\documentclass[12pt, Russian]{article}
\usepackage[russian]{babel}
\evensidemargin=\oddsidemargin \textwidth=165mm \textheight=238mm
\usepackage[dvips]{graphicx}
\usepackage{latexsym}
\usepackage{amssymb}
\begin{document}
\null
\begin{flushleft}\bf{\small}
\end{flushleft}
\vskip 5mm

\begin{center}
    \Large \bf Existence and smoothness of the solution to the Navier-Stokes
   equation
       \end{center}
\begin{center}
        {\sl Dr. Bazarbekov Argyngazy B.     \\}
        {\it Kazakh National State University, EKSU.~ e-mail : arg.baz50@mail.ru}

\end{center}
 \textbf{Abstract~.}\begin{scriptsize} A fundamental problem in
analysis is to decide whether a smooth  solution exists for the
Navier-Stokes equations in three dimensions . In this paper we shall
study this problem. The Navier-Stokes equations are given
by:$~~~~~~~~~~~$ $u_{it}(x,t) - \rho \triangle u_i(x,t) - u_j(x,t)~
u_{ix_j}(x,t) + p_{x_i}(x,t) = f_i(x,t)$ , $div~\textbf{u}(x,t) =
0~~ ,~i = 1,2,3$~ with initial conditions $\textbf{u}|_{(t=0)\bigcup
\partial\Omega} = 0$. We introduce the unknown vector-function:
$\big(w_i(x,t)\big)_{i=1,2,3}:~ u_{it}(x,t) - \rho \triangle
u_i(x,t)-\frac{d p(x,t)}{d x_i} = w_i(x,t)$~ with initial
conditions:~$u_i(x,0)= 0,$ $u_i(x,t)\mid_{\partial \Omega} = 0$. The
solution ~$u_i(x,t)$ ~of this problem is given by: $u_i(x,t) =
\int_0^t\int_{\Omega}G(x,t;\xi,\tau)~\Big(w_i(\xi,\tau) + \frac{d
p(\xi,\tau)}{d \xi_i}\Big) d \xi d \tau$~where $G(x,t;\xi,\tau)$ is
the Green function.~We consider the following Navier- Stokes -2
problem : find a solution~$\textbf{w}(x,t)\in \textbf{L}_2(Q_t),
p(x,t): p_{x_i}(x,t) \in L_2(Q_t)$~of the system of
equations:~~$w_i(x,t) - G \Big(w_j(x,t) + \frac{d p(x,t)}{d
x_j}\Big)\cdot G_{x_j}\Big(w_i(x,t) + \frac{d p(x,t)}{dx_i}\Big) =
f_i(x,t),i = 1,2,3$ satisfying  almost every where on $Q_t.$ ~Where
the vector-function $\big(p_{x_i}(x,t)\big)_{i=1,2,3}$ is defined by
the vector-function $\big(w_i(x,t)\big)_{i=1,2,3}$. ~~Using the
following estimates for the Green function: $\big|G(x,t;\xi
,\tau)\big| \leq c/ (t - \tau)^{\mu}\cdot |x - \xi|^{3-2\mu};
\big|G_{x}(x,t;\xi ,\tau)\big| \leq c/(t - \tau)^{\mu}\cdot|x -
\xi|^{3-(2\mu-1)} (1/2 < \mu < 1),$ from this system of equations we
obtain: $w(t) < f(t) + b \Big(\int_0^{t}\frac{w(\tau) + p(\tau)} {(t
- \tau)^{\mu}} d \tau\Big)^2$
 where $\mu: 5/8 < \mu < 1 , b  = const;~ w(\tau) = \|\textbf{w}(x,\tau)\|_{L_2(\Omega)};~f(t) =
\|\textbf{f}(x,t)\|_{L_2(\Omega)}, p(\tau) =
\sum_1^3\|\frac{\partial p(x,\tau}{\partial
x_i}\|_{L_2(\Omega)}.$~Using the estimate:~ $p(t)~<~c~w(t)$~ from
this inequality we infer: $w(t) < f(t)~ +~ b\cdot  \Big(\int_0^t
\frac{w(\tau)}{(t - \tau)^{\mu}}~d \tau\Big)^2$ where b is  real
number.~~ After the replacements of the functions  $\int_0^t
\frac{w(\tau) d \tau}{(t - \tau)^{\mu}}~=~w_1(t)$~ and ~ $z(t) =
z(0)e^{-k \int_0^t w_1(\tau) d \tau} $~ this inequality will accept
the following form:
$\frac{1}{k}~\int_0^t\frac{\frac{1}{z(\tau)}~\frac{d^2 z(\tau)}{d
\tau^2}}{(t - \tau)^{1- \mu}}~d \tau~+~f^2(t)~>~0$ where~$\mu: 5/8 <
\mu < 1$ is a real number. This is analogue of the replacement of
function by \textbf{Riccati} : $ z(t)=- \frac{1}{b}\cdot
\frac{u'(t)}{u(t)}$ for the solution of the following ordinary
nonlinear equation : $\frac{d z(t)}{d t} = f(t) + b z^2(t);~z(0) = 0
[10p.41].$ From the last inequality we obtain the a priori estimate:
$\|\textbf{w}(x,t)\|_{L_2(Q_t)} <
\sqrt{2}\|\textbf{f}(x,t)\|_{L_2(Q_t)}$ where $Q_t = \Omega\times
[0,t]$, t > 0 is an arbitrary real number. ~By the well known
Leray-Schauder's method and this a priori estimate  the existence
and uniqueness of the solution  $\textbf{u}(x,t): \textbf{u}(x,t)
\in \textbf{W}_{2}^{2,1}(Q_t)\bigcap \textbf{H}_2(Q_t)$ is proved.~
We used the nine known classical theorems.
\end{scriptsize}

                                               \label {Baizhuman}

\vskip 1mm

\begin{scriptsize} 2000 Mathematics Subject Classification . Primary :~$35~K~55$ ;
Secondary : $46 E 35$;  Keywords :~The Holder's inequalities,
Theorems: of Volterra V., of Hardy-Littlewood  , of Sobolev S.L., of
Lerau-Schauder, of Weyl H., of Abel-Carleman , of Riccati, Gronwall'
Lemma. Estimates for the Green function , Gamma function ,
projection operators , a priori estimate.\end{scriptsize}

\textbf{1.~Introduction}.~ The Navier-Stokes equations are given by

$$\frac{\partial u_i(x,t)}{\partial t}
~~-~~\rho~\triangle~u_i(x,t)~~+~\sum_{j=1}^{n}u_j(x,t)~
\frac{\partial u_i(x,t)}{\partial x_j}~+~\frac{\partial
p(x,t)}{\partial x_i}~~=~~f_i(x,t)$$

\[div~~u(x,t)~=~0~;~~i~=~1,...,n\]
with the initial condition :$$u (x,0) =~u_0(x)    \eqno(1.2)$$
where~$u (x,t)~=~(u_i(x,t))_{i =1,2,3}$ and $p(x,t)\in R$ are the
unknown velocity vector and pressure~defined for position ~~$x\in
R^3$~~and time ~~$t \geq 0$ . Here,  $u_0(x)$ is a given
divergence-free vector field on ~~$R^n$~,~ ~$f_i(x,t)$~~are the
components of a given  externally force ,~~$\rho > 0$ ~~is a
positive coefficient ~ and ~~$\triangle = \sum
_{i=1}^{3}\frac{d^2}{d x_i^2}$~~is the Laplacian in the space
variables .

Starting with Lerau [1] , important progress has been made in
understanding weak solutions of the Navier-Stokes equations.For
instance, if (1.1) and (1.2) hold, then for any smooth vector field
 $\varphi(x,t)=(\varphi_i(x,t))_{i =1,2,3}$, compactly supported
in~$R^n \times (0,\infty)$,  a formal integration by parts yields
$$\int\int_{R^n\times R} u(x,t)
\frac{\partial\varphi}{\partial t}~-~\sum_{i,j}\int\int_{R^n\times
R} u_i(x,t)u_j(x,t)\frac{\partial\varphi_i}{\partial x_j} dx dt =$$
$$=-\rho\int\int_{R^n\times R} u(x,t) \bigtriangleup \varphi
dx dt+\int\int_{R^3\times R}f(x,t)\varphi dx dt
 - \int\int_{R^n\times R}p(x,t) (div \varphi) dx dt \eqno(1.3)$$

Note that (1.3) makes sense  for~$u \in L_2 , f \in L_1 , p \in L_1
$~ whereas (1.1) makes sense only if~ u(x,t)~is twice differentiable
in~x~.~Similarly, if ~~$\varphi(x,t)$~is a smooth function,
compactly supported in ~~$R^n\times (0,\infty)$,~~then a formal
integration by parts and (1.2) imply: $$\int\int _{R^n\times R}
u(x,t)\bigtriangledown_{x}\varphi(x,t)~dx dt  ~~=0   \eqno(1,4)$$. A
solution (1.3),(1.4) is called a weak solution of the Navier-Stokes
equations. Leray in [1] showed that  the Navier-Stokes equation
(1.1), (1.2), (1.3) in three space dimensions always have a weak
solution~~(u(x,t) , p(x,t)). The uniqueness of weak solutions of the
Navier-Stokes equation is not known.~ In two dimensions the
existence, uniqueness and smoothness of weak solutions have been
known for a long time (R.Temam [2],O. Ladyzhenskaya [3],I.
Lions[4]).

In three dimensions, this questions studied for  the initial
velocity ~~$u_0(x)$~ satisfying  a smallness condition . For the
initial data~~$u_0(x)$~~not assumed to be small , it is known that
the existence of smooth weak solutions holds if the time interval
~~$[0,\infty)$~ is replaced by a small time interval ~~$[0 , T)$
depending on the initial data .

A fundamental problem in analysis is to decide whether a smooth
solution exists for the Navier-Stokes equations in three dimensions.

 \vskip 1mm \textbf{2.~~Results}

Let~$\Omega \subset R^3$~be a finite domain bounded by the Lipchitz
surface ~$\eth \Omega.$~$Q_{t} = \Omega \times[0,t] , ~x
=(x_1,x_2,x_3)$~and~$\textbf{u}(x,t)~=~ (u_i(x,t)_{i = 1,2,3}~,$
$\textbf{f}(x,t) = (f_i(x,t)_{i =1,2,3}$ are  vector-functions.~
Here~t > 0 is an arbitrary real number.~ The Navier-Stokes equations
are given by
$$\frac{\partial u_i(x,t)}{\partial t}
~-~\rho~\triangle~u_i(x,t)~-~\sum_{j=1}^{3}u_j(x,t)~ \frac{\partial
u_i(x,t)}{\partial x_j}~+~\frac{\partial p(x,t)}{\partial
x_i}~~=~f_i(x,t)~\eqno(2.1),$$
$$div~\textbf{u}(x,t)~=~\sum_{i=1}^3\frac{\partial u_i(x,t)}{\partial x_i}~=~
0~~,i~=~1,2,3$$\textbf{The Navier-Stokes problem~1.} Find a
vector-function  $\textbf{u}(x,t)~=~(u_i(x,t))_{i=1,2,3} : \Omega
\times [0,t] \rightarrow R^3$ and a scalar function $p(x,t):
 \Omega \times [0,t] \rightarrow R^1$  satisfying  the equation
(2.1) and   the following initial condition
$$\textbf{u} (x,0) =~0~~, ~~~ \textbf{u}(x,t)\mid_{\partial \Omega \times [0,t]}~=~0
\eqno(2.2)$$Let $p > 1, r > 1$ be real numbers. We shall use the
following functional spaces.
$~~~~~~~~~~~~~~~~~~~~~~~~~~~~~~~~~~~~~~~~~~~~~~~~~~~~~~~~~~~~~~~~~~~~~~~~~~~~~~~~~~~~~~$
$L_{p,r}(Q_t)$ ~~ is the Banach space with the norm~ [3 p.33]
$$\|u(x,t)\|_{L_{p,r}(Q_{t})}~~=~~\Big[ \int_0^t
\Big(\int_{\Omega}|u(x,t)|^{p} dx \Big)^{r/p}~dt
\Big]^{1/r}~,~L_{p,p}(Q_t)~=~~L_p(Q_t).$$ $W^{2,1}_p (Q_t)$ ~ is the
Banach space supplied by the norm
$$\|u(x,t)\|_{W^{2,1}_{p}(Q_t)} = \Big[\|u\|^p_{L_{p}(Q_t)} + \|u_t\|^p_{L_{p}(Q_t)} +
\|u_x\|^p_{L_{p}(Q_t)} + \|u_{xx}\|^p_{L_{p}(Q_t)}\Big]^{1/p}$$
$\textbf{L}_{p,r}(Q_t)$ ~~~is the Banach vector-space with the norm
$$\|\textbf{u}(x,t)\|_{L_{p,r}(Q_t)}~~=~~\sum_{i=1}^3\|u_i(x,t)\|_
{L_{p,r}(Q_t)}$$
$\textbf{L}_{2}(Q_t)$ ~~is the Hilbert
vector-space with the inner product
$$(\textbf{u}(x,t) , \textbf{v}(x,t))_{L_{2}(Q_t)}~=~
\sum_{i=1}^3(u_i(x,t),v_i(x,t))_{L_{2}(Q_t)}$$
$\textbf{V}_0(Q_t)~\big(\textbf{V}(Q_t)\big)$ are the vector-spaces
of  smooth functions
 $$V_0(Q_t) = \{\textbf{u}(x,t)\in C^2 (\overline{Q_t}),~div~ \textbf{u}(x,t) = 0~,
 \textbf{u}\cdot \textbf{n}|_{\partial \Omega} = \sum_{i=1}^3 u_i(x,t) cos (\textbf{n},x_i)|_{\partial \Omega} = 0\},$$
$$V(Q_t) = \big(\{\textbf{u}(x,t)\in C^2 (\overline{Q_t}): div~
\textbf{u}(x,t) = 0\}\big)$$$\textbf{H}_2(Q_t)$ is the closure of ~
$\textbf{V}_0(Q_t)$~ in the norm of $\textbf{L}_2(Q_t).$  [2
p.13]~I.e.$$\textbf{H}_2(Q_t) = \{\textbf{u}(x,t) :
\textbf{u}(x,t)\in \textbf{L}_2(Q_T) , div~\textbf{u}(x,t) = 0,
\textbf{u}\cdot \textbf{n}|_{\partial\Omega}  = 0\}$$
$\textbf{E}_2(Q_t)$ is the closure of ~ $\textbf{V}(Q_t)$~ in the
norm of $\textbf{L}_2(Q_t).$~~[2 p. 13] I.e.
$$\textbf{E}_2(Q_t) = \{\textbf{u}(x,t) : \textbf{u}(x,t)\in
\textbf{L}_2(Q_t) , div~\textbf{u}(x,t) = 0\}$$  It is obvious that
$\textbf{H}_2(Q_t) \subseteq \textbf{E}_2(Q_t).$~ Further, we shall
denote  the vector-functions and  vector-spaces by bold type.
~~~~~~The following is  principal result.

\textbf{Theorem~2.1.}~~For any right-hand side ~$\textbf{f}(x,t) \in
\textbf{L}_2(Q_t)$~~in equation (2.1) and for any real numbers
~~$\rho > 0 , t > 0,$~~the Navier-Stokes problem-1 has a unique
smooth solution ~$\textbf{u}(x,t) : \textbf{u}(x,t)
\in\textbf{W}_2^{2,1}(Q_t) \cap \textbf{H}_2(Q_t)$ and a scalar
function  $p(x,t): p_{x_i}(x,t) \in L_2(Q_t)$ satisfying  (2.1)
almost everywhere  on $Q_t,$ and  the following estimates are valid:
$$\|\textbf{u}(x,t)\|_{W^{2,1}_2(Q_t)}~\leq~c~\|\textbf{f}\|_{
L_2(Q_t)}~,~\Big\|\frac{\partial p(x,t)}{\partial x_i}
\Big\|_{L_2(Q_t)} ~~\leq c~\|\textbf{f}\|_{L_2(Q_t)} \eqno (2.3)$$
Here and bellow by symbol~$c,$~~we denote a generic constant ,
independent on the solution and  right-hand side  whose value is
inessential to our aims, and it may change from line to line.

\textbf{Remark~2.1.}~~The case when the right-hand
side~$\textbf{f}(x , t)$~has a small norm or a time~ $t \ll 1$~~
$(or~\rho \gg 1)$~~~is well-known and so not interesting. But in
Theorem 1 ~$\textbf{f}(x,t) \in \textbf{L}_2(Q_t)$~ is an arbitrary
vector-function and ~t > 0~,~$\rho > 0$~ are arbitrary real numbers
.  In recent paper [6] Ladyzhenskaja formulates the Navier-Stokes
problem as in the formulas (2.1) - (2.2) and in Theorem 1 . For
simplicity, we consider the Navier-Stokes problem for the
homogeneous case (i.e. $u (x,0) =~0~,~ u(x,t)\mid_{\partial
\Omega}~=~0.$).  We consider the inhomogeneous case (i.e. $u (x,0)
=~u_0(x)~,~ u(x,t)\mid_{\partial \Omega}~=~0$) in Section 4.

\textbf{Definition~2.1.} A vector-function $\textbf{u}(x,t):
\textbf{u}(x,t)|_{(t=0)\cup \partial \Omega} = 0$ and a scalar
function p(x,t) are called a \textbf{smooth} solution to the
Navier-Stokes problem-1, if $\textbf{u}(x,t)
\in\textbf{W}_2^{2,1}(Q_t) \cap \textbf{H}_2(Q_t)$ and $p_{x_i}(x,t)
\in L_2(Q_t).$

We adduce the well-known definition of the Hopf solution to the
Navier-Stokes equation.

\textbf{Definition 2.2} (the Hopf's solution). ~ Let a right-hand
side ~$\textbf{f}(x,t) \in \textbf{L}_2(Q_t).$ A vector-function
$\textbf{u}(x,t) \in \textbf{L}_2([0,t];H_0(\Omega))\cap
\textbf{L}^{\infty}([0,t]; L_2(\Omega))$ is called the Hopf's
solution, if the following equality [2 p.225].$$\frac{\partial
\big(\textbf{u}(x,t) , \textbf{v}(x)\big)}{\partial
t}~+~\rho~\big(\textbf{u}_x(x,t) , \textbf{v}_x(x)\big)
~-~\sum_{i=1}^{3}u_j(x,t) \big(\textbf{u}_{x_j}(x,t)~ ,
~\textbf{v}(x)\big) = \int_{\Omega} \textbf{f}(x,t) \cdot
\textbf{v}(x) d x$$ is fulfilled for all vector-functions
$\textbf{v}(x)\in\textbf{H}_0^1(\Omega)=\{\textbf{u}(x):
div~\textbf{u}(x)=0,\textbf{u}(x)|_{\partial \Omega}=0 ,
\textbf{u}(x),\textbf{u}_{x_i}(x) \in \textbf{L}_2(\Omega)[2
p.24]\}.$

\textbf{Remark 2.2.} By Theorem 2.1 it follows that the Hopf's
solution is the smooth solution.$\blacktriangleleft$

\textbf{For the proof of Theorem 2.1 we shall use the following
known propositions.}

\textbf{Theorem of Weyl H.} In the book [2 p.22] the following
equalities are proved:

$\textbf{L}_2(Q_t) = \textbf{H}_2(Q_t) \oplus \textbf{G}_2(Q_t)$ ~
where $\textbf{H}_2(Q_t) = \{\textbf{u}(x,t) : \textbf{u}(x,t)\in
\textbf{L}_2(Q_t) , div~\textbf{u}(x,t) = 0, \textbf{u}\cdot
\textbf{n}|_{\partial \Omega \times [0,t]} = 0 \}.$
$\textbf{G}_2(Q_t) = \{\textbf{u}(x,t) : \textbf{u}(x,t) \in
\textbf{L}_2(Q_t), \textbf{u}(x,t) = \textbf{grad} p(x,t):
p_{x_i}(x,t) \in L_2(Q_t)\}.$ ~~~ I.e. for any $\textbf{f}(x,t) \in
\textbf{L}_2(Q_t),$ the following equality:
$\textbf{f}(x,t)~=~H\big(\textbf{f}(x,t)\big) +~G
\big(\textbf{f}(x,t)\big)$ is valid ~ where $H :
\textbf{L}_2(Q_t)\Rightarrow \textbf{H}_2(Q_t)~ ,~G:
\textbf{L}_2(Q_t)\Rightarrow \textbf{G}_2(Q_t)$- are the projection
operators.

\textbf{Proposition~1.} ~(The Holder inequality).~ ~Let ~$p_1 > 1 ,
p_2 > 1 ; r_1 > 1~,~r_2 > 1$~be a real numbers. ~Then, the following
Holder inequality is valid ~~~~~~~~~[5~p.75].
$$\|u(x,t) v(x,t)\|_{L_{\frac{p_{1} p_{2}}{p_1 + p_2},\frac{r_1 r_2}{r_1
+r_2}}(Q_t)}~~~\leq~\|u(x,t)\|_{L_{p_1,r_1}(Q_t)}~
\|v(x,t)\|_{L_{p_2 ,r_2}(Q_t)}~\eqno(2.4)~$$ \textbf{Proposition
~2.} ~(The system equations  of  Volterra V.)~~ On the space of
vector-functions ~ $\textbf{u}(x,t) = (u_i(x,t))_{i=1,2,3} \in
$\textbf{L}$_2(Q_t)$~we shall consider the following system of
nonlinear integral equations of Volterra: ~ [7 p.59 , p.62]
$$u_l(x,t)~-~\sum_{s=1}^{3}\int_0^t \int_{\Omega}
K_{l,s}\big(x,t;\xi,\tau
;\textbf{u}(\xi,\tau)\big)~u_s(\xi,\tau)~d\xi~d\tau
~~=~f_l(x,t)\eqno(2,5)$$ l~=~1 ,2 ,3 .~~~~~Or, in the vector form
$$\textbf{u}(x,t)~~-~~K~\textbf{u}(x,t)~~~=~~~\textbf{f}~(x,t)~\in~
\textbf{L}_2(Q_T)~\eqno(2,6)$$ This system of equations under some
conditions to the nonlinear kernel
$K(x,t;\xi,\tau;\textbf{u}(\xi,\tau))$  has been studied   in the
book [7 p.61]. We shall study this nonlinear system of equations by
using the theorem of Leray J., Schauder J.

\textbf{Proposition~3.}~(Theorem of Hardy G.H., Littlewood J.E.) Let
$\mu :0 < \mu < 1$ be a real number. We shall consider the following
operator of the fractional integration~~~
$J^{\mu}u(t)~=~\int_0^t~\frac{u(\tau)~d \tau}{(t - \tau)^{\mu}}$.
Then:

\textbf{a)}~If~~$1~<~p~<~\frac{1}{1-\mu},$~~then the
operator~~$J^{\mu}$~~is bounded from the space ~~$L_p(0,t)$~into the
space~ $L_{q}(0,t)$~~where~~$q~=~\frac{p}{1-p\cdot
(1-\mu)}$~~and~~$\|J^{\mu}u(t)\|_{L_q(0,t)}~\leq~c~\|u(t)\|_{L_p(0,t)}.$
[8 p.64].

\textbf{Proposition~4.}~(Theorem of Sobolev S.L.)~~Let a
function~u(x)~ be  represented as the potential of a function f(x)
,~i.e. ~~~$u(x)~~=~~\int_{\Omega}~~\frac{f(\xi)~d~\xi}{|x~-~\xi|^{3
- \lambda}}~~~~~~\lambda~>~0.$~~~~[3~p.32].~~~~~~~ Then

\textbf{a)}~If~~$0~<~\lambda~<~3/p$~~and~$f(x)~\in~L_p(\Omega),$~
then~$u(x)~\in~L_q(\Omega)$~where~~ $q~\leq~\frac{3p}{3 -
p~\lambda}$~ and
 $\|u(x)\|_{L_q(\Omega)}~\leq~c~\|f(x)\|_{L_p(\Omega)}$.

\textbf{b)}~If~$~\lambda =~3/p$~and~$f(x)~\in~L_p(\Omega)$,~
then~$u(x)\in L_{\infty}(\Omega)$ and~$\|u(x)\|_{L_{\infty}(\Omega)}
\leq~c~\|f(x)\|_{L_p(\Omega)}.$

\textbf{Proposition~5.}~We shall consider the following problem on
the domain $Q_t$ :
$$ u_t(x,t)~~-~~\rho~\triangle~u(x,t)~=~g(x,t),~~~~~~
u (x,0) =~0~~, ~ u(x,t)\mid_{\partial \Omega\times [0,t]}~=~0
\eqno(2,7)$$ The solution ~$u(x,t)$ ~of this problem is given by
$$u(x,t)~=~G g(\xi,\tau)~=~\int_0^t\int_{\Omega}G(x,t;\xi,\tau)~g(\xi,\tau)~
d \xi~d \tau   \eqno(2,8)$$ where~~$G(x,t;\xi,\tau)$~ is the Green
function~for $Q_t$~.~~The construction of the Green function is
resulted in book [9~p.111]. The following estimates are valid.[9~
p.170]
$$\big|~G (x,t;\xi ,\tau)\big|~~\leq~~\frac{const}{(t -
\tau)^{\mu}}~\frac{1}{|x - \xi|^{3-2 \mu}}~~,~~~~0~<~\mu~<~1
\eqno(2,9)$$
$$\big|~\frac{\partial}{\partial~x}~G (x,t;\xi ,\tau)\big|~~\leq~~\frac{const}{(t -
\tau)^{\mu}}~\frac{1}{|x - \xi|^{3-(2 \mu~-~1)}}~~,~~~~1/2~<~\mu~<~1
\eqno(2,10)$$ From estimates in Propositions~3~,~4~and~the estimates
(2,9),(2,10)~ follow that:

\textbf{a)}~If~$g(x,t)~\in~L_2(Q_t)$~and~$\mu = 5/8,$~~then
$$G g(x,t)~\in~L_{12,8}(Q_t)~~,~~G_x
g(x,t)~=~\frac{\partial G g(x,t)}{\partial x}~\in~L_{12/5,8}(Q_t),
\eqno(2,11)$$
$$\big\|G g(x,t)~G_x g(x,\tau)\big\|_{L_{2,4}(Q_t)}~\leq~
\big\|G g(x,t)\big\|_{L_{12,8}(Q_t)}~\big\|G_x
g(x,\tau)\big\|_{L_{12/5,8}(Q_t)}~\leq~ ~c~\|g(x,t)\|_{L_{2}(Q_t)}^2
$$ \textbf{b)}~~~If
~$g_1(x,t),~g_2(x,t)~\in~L_2(Q_t)$~and~$g_i(t)~=~\|g_i(x,t\|_{L_2(\Omega)}~i=
1, 2$,~then for any ~$\mu: 5/8 < \mu < 1$~the following estimates
are valid:
$$\Big\|G g_i (x,t)\Big\|_{L_{\frac{6}{3-4\mu}}(\Omega)} \leq c ~
\int_0^t\frac{g_i(\tau)~d \tau}{(t - \tau)^{\mu}};~~~~~~;\Big\|G_x
g_i (x,t)\Big\|_{L_{\frac{6}{5-4\mu}}(\Omega)} \leq c ~
\int_0^t\frac{g_i(\tau)~d \tau}{(t - \tau)^{\mu}} $$
$$\Big\|G g_1(x,t)~G_x g_2(x,\tau)\Big\|_{L_{\frac{3}{4(1-\mu)}}(\Omega)}~\leq
~~ c~\int_0^t\frac{g_1(\tau)~d \tau}{(t - \tau)^{\mu}}\cdot
\int_0^t\frac{g_2(\tau)~d \tau}{(t - \tau)^{\mu}}   \eqno(2.12.1)$$
and for ~ $\mu:~\frac{5}{8} \leq \mu < 1$~  follows that~ $2 \leq
\frac{3}{4(1-\mu)}.$ ~Therefore,
$$\Big\|G g_1(x,t)~G_x
g_2(x,\tau)\Big\|_{L_{2}(\Omega)}~\leq$$$$\leq~c~ \Big\|G
g_1(x,t)~G_x g_2(x,\tau)\Big\|_{L_{\frac{3}{4(1-\mu)}}(\Omega)}~\leq
~~ c~\int_0^t\frac{g_1(\tau)~d \tau}{(t - \tau)^{\mu}}\cdot
\int_0^t\frac{g_2(\tau)~d \tau}{(t - \tau)^{\mu}} \eqno(2.12.2)$$

$\blacktriangleright$~\textbf{a)}~We shall  prove the inequality
(2.11).~From Holder inequality (2,4) (with~$p=2, p_1 = 12 , r_1 =
8)$ and $(p=2, p_2 = 12/5 , r_2 = 8)$ ,~  the first inequality of
(2.11) follows. From the  estimate (2.9) with ~ $(\mu = 5/8)$~we
have
$$\big\|G g(x,\tau)\big\|_{L_{12,8}(Q_t)}~=~\Big\|\int_0^{\tau}\int_{\Omega}
G(x,\tau;\xi,\tau_1)~ g(\xi,\tau_1)~ d \xi~d
\tau_1\Big\|_{L_{12,8}(Q_t)}~ \leq $$$$\leq~c~
\Big\|\int_0^{\tau}\frac{1}{(\tau -
\tau_1)^{5/8}}\Big\|\int_{\Omega}\frac{g(\xi , \tau_1) d \xi}{|x -
\xi|^{3 -~5/4}}\Big\|_{L_{12}(\Omega)}d
\tau_1\Big\|_{L_8(0_t)}~\leq~$$ $$\leq c~\Big\|\int_0^{\tau}
\frac{\|g(x,\tau_1)\|_{L_2(\Omega)} d \tau_1 } {(\tau-
\tau_1)^{5/8}}\Big\|_{L_8(0_t)}~\leq~c~ \|g(x,t)\|_{L_2(Q_t)}$$ Here
we used the fact that from Proposition $4_a$~\big(with $p=2,
\lambda=5/4; q = \frac{6}{3-2 \lambda} = 12$\big)   follows the
inequality $\Big\|\int_{\Omega}\frac{g(\xi , \tau_1) d \xi}{|x -
\xi|^{3 -~5/4}}\Big\|_{L_{12}(\Omega)} \leq~c
\|g(x,\tau_1)\|_{L_2(\Omega)}$~ and
since~$\|g(x,\tau_1)\|_{L_2(\Omega)} \in L_2(0,t),$~then~from
Proposition $3$~(with $p=2;\mu = 5/8;~q = \frac{p}{1-p~(1-\mu)} =
8$) the following inequality
$\Big\|\int_0^{\tau}\frac{\|g(x,\tau_1)\|_{L_2(\Omega)} d \tau_1 }
{(\tau-\tau_1)^{5/8}}\Big\|_{L_8(0,t)}~\leq~c~\|g(x,t)\|_{L_2(Q_t)}$~follows.

Using the estimate  (2,10) (with $\mu: \frac{5}{8} \leq \mu < 1$) ,
the Proposition $4_a$~\big(with $\lambda=1/4; q = \frac{6}{3- 2
\lambda} = \frac{12}{5}$\big) and  Proposition $3$~(with $p=2;\mu =
5/8;~q = \frac{p}{1-p~(1-\mu)} = 8$)   the following estimate
~$\big\|\frac{\partial G g(x,t)}{\partial x}\big\|_{L_{12/5,8}(Q_t)}
\leq c \|g(x,t)\|_{L_{2}(Q_t)}$ is proved similarly.~~ By these
estimates  the second estimate of (2.11) follows.~ The inequality
(2.11) is proved.

\textbf{b)}~The proofs of the inequalities in (2.12.1), (2.12.2)~
follow from the estimates of the Green function (2,9), (2,10),
Proposition $4_a$ and is similar to the proof of \textbf{a}). The
parameters in the second inequality of (2,12.1) satisfy  all
conditions of parameters  $p_1 , p_2$ ~in (2,4) of Proposition 1.
$\blacktriangleleft$

\textbf{Proposition 6.} ( Theorem of Leray J., Schauder J.)~Let  А
be a compact  nonlinear operator on $L_2(Q_t)$.~~If  every possible
solution to the following equation  [3 p.42]
$$    w(x,t) ~~+~~A w(x,t)~~= ~~f(x,t)$$do not fall outside the bounds
of some sphere ~$|w(x,t)|_{L_2(Q_t)} \leq c,$~ then  for any
right-hand side $f(x,t)\in L_2(Q_t)$ the equation has at least one
solution in this sphere.

\textbf{Proposition 7.}(The equation of Abel-Carleman.)~The equation
of Abel-Carleman is set by the following formulas [8 p.39]:
$$\int_0^t \frac{u(\tau)}{(t-\tau)^{\mu}}~d \tau = f(t);~~~~~~~~~~
u(t) = \frac{sin \pi \mu}{\pi} \frac{d}{d t} \int_0^t
\frac{f(\tau)}{(t- \tau)^{1-\mu}} d \tau \eqno(2.13)$$ Let ~~~$-
\infty < \mu_1 < 1 ;~~- \infty < \mu_2 < 1$.~~~~Then
$$\int_0^t \frac{d \tau}{(t-\tau)^{\mu_1}} d \tau \int^{\tau}_0 \frac{g(\tau_1) d
\tau_1}{(\tau-\tau_1)^{\mu_2}}~  = \int_0^t g(\tau_1)~d \tau_1
\int_{\tau_1}^t \frac{d
\tau}{(t-\tau)^{\mu_1}~(\tau-\tau_1)^{\mu_2}}  = $$
$$= \Gamma_{\mu_1}^{\mu_2}\int_0^t
\frac{g(\tau)~d \tau}{(t-\tau)^{\mu_1+\mu_2 - 1}}~;~~~~~~~~~~~~
;\Gamma_{\mu_1}^{\mu_2}~=~\frac{\Gamma(1-\mu_1)
\Gamma(1-\mu_2)}{\Gamma(2-\mu_1-\mu_2)}
   \eqno(2.14)$$

\textbf{Proposition ~8.} (The linear Navier-Stokes equation) Let t
> 0 be an arbitrary real number. ~We consider the
following linear Navier-Stokes problem on the domain $Q_t$:~~[3,
p.95]

$$\frac{\partial u_i(x,t)}{\partial t}
~-~\rho~\triangle~u_i(x,t)~-~\frac{\partial p(x,t)}{\partial
x_i}~=~w_i(x,t)~\eqno(2.15),$$
$$div~\textbf{u}(x,t)~=~0~;~~\textbf{u} (x,0) =~0~,~
\textbf{u}(x,t)\mid_{\partial \Omega \times [0,t]}~=~0$$ For this
problem in the manuscript of author [12] is received the explicit
expression to the  pressure function  p(x,t) , depending  on the
right-hand side $w_i(x,t)$:
$$p(x,t) ~ =~ - T\ast\bigtriangleup^{-1}\ast\int_0^t \int_{\Omega}\sum_1^3
\frac{dG(x,t;\xi,\tau)}{d x_i}w_i(\xi,\tau) d \xi d \tau ~=
\eqno(2,16)$$
$$=-\frac{d }{d t}
\triangle^{-1}\ast \int_0^t \int_{\Omega}\sum_1^3
\frac{dG(x,t;\xi,\tau)}{d x_i}w_i(\xi,\tau) d \xi d \tau + \rho
\cdot \int_0^t \int_{\Omega}\sum_1^3 \frac{dG(x,t;\xi,\tau)}{d
x_i}w_i(\xi,\tau) d \xi d \tau$$ where~ $T\ast u(x,t) = \frac{d}{d
t}u(x,t) - \rho \triangle_x u(x,t)$ is the parabolic operator ~,~
$\triangle^{-1}$ is the inverse  operator to  Dirichlet problem for
Laplase equation on the domain  $\Omega$. ~$G(x,t;\xi,\tau)$ is the
Green function of Dirichlet problem  for the parabolic equation on
the domain $Q_t = \Omega\times [0,t]$ [9 p.106].~If $\textbf{w}(x,t)
\in \textbf{L}_2(Q_t)$, then $\int_0^t\int_{\Omega} G(x,t;\xi,\tau)
\textbf{w}(\xi,\tau) d \xi d \tau ~\in~\textbf{W}^{2,1}_2(Q_t)$.~~
It is obvious that:
$$\Big|\int_0^t \int_{\Omega}\sum_1^3 \frac{dG(x,t;\xi,\tau)}{d
x_i}w_i(\xi,\tau) d \xi d \tau \Big|_{W^{1,1}_2(Q_t)}~<~~c\cdot
|\textbf{w}(x,t)|_{L_2(Q_t)}$$From this estimate follows the
following estimate:
$$\int_0^t\Big(\Big|\frac{d}{d t} \bigtriangleup^{-1}\ast\int_0^{\tau}
\int_{\Omega}\sum_1^3 \frac{d G(x,t;\xi,\tau_1)}{d
x_i}w_i(\xi,\tau_1) d \xi d \tau_1\Big|_{W_2^1(\Omega)}\Big)^2 d
\tau ~<~c~|\textbf{w}(x,\tau)|^2_{L_2(Q_t)}$$ We consider the
formula (2,16) in detail. ~ It is known that the following classical
problem for the parabolic equations:~$u_t - \rho \triangle u(x,t) =
f(x,t) \in L_2(Q_t), ~u_{t=0\cup \partial \Omega} = 0$~has the
unique solution $u(x,t) \in W_2^1(Q_t)$ ~and $\|u(x,t\|_{W_2^1(Q_t)}
< c \|f(x,t)\|_{L_2(Q_t)}$ where $Q_t = \Omega \times
[0,t].$~~Therefore, from the formula (2,16)  follow the following
estimates:
$$\Big\| \frac{d}{d x_i} \frac{d }{d t} \triangle^{-1}\ast \int_0^t \int_{\Omega}\sum_1^3
\frac{dG(x,t;\xi,\tau)}{d x_i}w_i(\xi,\tau) d \xi d
\tau\Big\|_{L_2(Q_t)} ~<~c \|w(x,t)\|_{L_2(Q_t)}$$
$$ \Big\|\frac{d}{d x_i} \rho \cdot \int_0^t \int_{\Omega}\sum_1^3
\frac{dG(x,t;\xi,\tau)}{d x_i}w_i(\xi,\tau) d \xi d
\tau\Big\|_{L_2(Q_t)} ~<~c \|w(x,t)\|_{L_2(Q_t)}$$

Differentiating the function of pressure  p(x,t) in formula  (2,16)
by  $x_i$, we find $\frac{d p(x,t)}{d x_i}, i = 1,2,3,$ depending on
the right-hand side $w_i(x,t)$. And by these estimates, integrating
on the domain $Q_t$, we obtain the following estimate:
$$\int_0^t\sum_1^3\Big|\frac{\partial p(x,\tau)}{\partial x_i}\Big|_{L_2(\Omega)}^2 d \tau~
 <~ c \cdot \int_0^t\sum_1^3|w_i(x,\tau)|_{L_2(\Omega)}^2 d \tau
\eqno(2,17)$$ By the formula (2,16) and estimate  (2,17) on the
vector space  $\textbf{w}(x,t) ~=~\Big(w_i(x,t)\Big)_{i = 1,2,3}
~\in~\textbf{L}_2(Q_t)$~we define the following linear and bounded
operator:
$$ P \Big(w_i(x,t)\Big)_{i = 1,2,3}~=~\Big(\frac{d p(x,t)}{d
x_i}\Big)_{i = 1,2,3}~\in ~\textbf{L}_2(x,t) \eqno(2,18)$$ where the
functions  $\frac{d p(x,t)}{d x_i}$~ is defined by the functions
$w_i(x,t)$ ~from the formula (2,16). $\blacktriangleleft$~~

Using the Green function $G(x,t;\xi,\tau)$, from the equation (2,15)
we find :
$$u_i(x,t)~=~\int_0^t\int_{\Omega}G(x,t;\xi,\tau)~\Big(w_i (\xi,\tau) + \frac{\partial p(\xi,\tau)}
{\partial \xi_i}\Big) d \xi d \tau \eqno(2,19)    $$ And present the
nonlinear Navier-Stokes equations (2,1) as:
$$w_i(x,t)~~-~\sum_{j=1}^{3}G \Big(w_j(x,t) + \frac{\partial p(x,t)}{\partial x_j}\Big)\cdot
 G_{x_j}\Big(w_i(x,t) + \frac{\partial p(x,t)}{\partial x_i}\Big)~~=~~f_i(x,t)~\eqno(2.20)$$
where~$G_{x_j}w_i(x,t) = \frac{\partial G w_i(x,t)}{\partial x_j}.$
~ Differentiating the pressure function  p(x,t) in the formula
(2,16) by $x_i$, we find  $\frac{d p(x,t)}{d x_i}, i = 1,2,3,$
depending on the functions $w_i(x,t)$, and \textbf{ substitute}
these functions to the equation (2,20). \textbf{Then, for the
definition of the three unknown functions
~$\big(w_i(x,t)\big)_{i=1,2,3}$ we obtain the three system of
nonlinear equations of Volterra (2,20)}. ~

\textbf{Remark 2,3.}~It is obvious that the vector function
$\Big(\frac{d p(x,t)}{d x_i}\Big)_ {i = 1,2,3}$ depends on the
vector function $(w_i(x,t))_{i=1,2,3}$~linearly. But we shall not
write the expressions of these depends,~we shall use the estimate
(2,17).$\blacktriangleleft$

\textbf{Navier-Stokes problem 2.} Find the vector-function
 $\big(w_i(x,t)\big)_{i=1,2,3} \in \textbf{L}_2(Q_t) $, satisfying
the equation (2,20) almost every where on $Q_t$.$\blacktriangleleft$

We will find the unknown vector-function $\textbf{w}(x,t)~=
~\big(w_i(x,t)\big)_{i=1,2,3} \in \textbf{L}_2(Q_t).$

\textbf{Theorem 2.2}~For any right-side $\textbf{f}(x,t) \in
\textbf{L}_2(Q_T)$ in the system of equations (2,20) there exists a
unique vector-function~$\textbf{w}(x,t) \in \textbf{L}_2(Q_T)$, ~
satisfying almost everywhere on $Q_T$,~ the system equations (2,20)
.~~~~And for any possible solution $\textbf{w}(x,t):
\|\textbf{w}(x,t)\|_{L_2(Q_T)} = \Big\|
\|\textbf{w}(x,t)\|_{L_2(\Omega)}\Big\|_{L_2(0,T)} =
\|w(t)\|_{L_2(0,T)} < \infty$~ to the basis equation (2,20), the
following a priori estimate is valid
$$\|\textbf{w}(x,t)\|_{L_2(Q_T)}~~~<~~~~\sqrt{2}\cdot  \|\textbf{f}(x,t)\|_{L_2(Q_{T})}   \eqno(2.21)$$
where  $Q_T = \Omega\times [0,T]$, ;~~ T > 0 is an arbitrary real
number.

\vskip 2mm \textbf{3.~Proof of Theorem 2.2.}

\textbf{The following is the key Lemma for the proof of Theorem
2.2.}

\textbf{Lemma~3.1.}~Let the vector-function~~$\textbf{g}(x,t)~=~(g_i
(x,t))_{i=1,2,3} \in~\textbf{L}_2(Q_t)$.~~~ I.e.
$|\textbf{g}(x,\tau)|_{L_2(Q_t)}~=~\sum_{i=1}^{i=3}(\int_{Q_t}g_i^2(x,\tau)~dx
d \tau)^{1/2}~<~\infty.$~~ We define the following vector-function:
$$G g_j(x,t) G_{x_j}\textbf{g}(x,t) = \Big(\sum_{j=1}^3
G g_j(x,t) G_{x_j} g_i(x,t)\Big)_{i=1,2,3} = \Big(\sum_{j=1}^3 G
g_j(x,t) \frac{\partial G g_i(x,t)}{\partial x_j}\Big)_{i=1,2,3}$$
Then for any ~t > 0 ~there exists~a constant~b > 0 independent on
\textbf{g}(x,t)~such that the following inequality is valid:
$$\Big\|G g_j(x,t) G_{x_j}\textbf{g}(x,t)\Big\|_{L_2(\Omega)}~ \leq~ b
\Big(\int_0^{t}\frac{g(\tau) d \tau}{(t - \tau)^{\mu}}\Big)^2
\eqno(3.1)$$ where~~$\mu:~\frac{5}{8} \leq\mu
<1$~,~$g(\tau)~=~\|\textbf{g}(x,\tau)\|_{L_2(\Omega)} =
\sum_{i=1}^{i=3}(\int_{\Omega}g_i^2(x,\tau)~dx)^{1/2}.$

$\blacktriangleright$~Let ~$\mu : \frac{5}{8} \leq \mu < 1.$~ By the
estimate (2.12.2) in Proposition 5 ~and the following
 inequality:~$\sum_{i,j=1}^3 c_i ~c_j\cdot a_i \cdot a_j \leq b~(\sum_{i=1}^3
a_i)^2$~~we have
$$\big\|G g_j(x,t) G_{x_j}\textbf{g}(x,t)\|_{L_2(\Omega)}~ \leq~
\sum_{i,j = 1}^3 \|G g_j(x,t) G_{x_j}g_i(x,t)\|_{L_2(\Omega)} ~<~$$
$$<~\sum_{i,j = 1}^3 c_i~c_j \int_0^t\frac{g_i(\tau)~d \tau}{(t - \tau)^{\mu}}\cdot
\int_0^t\frac{g_j(\tau)~d \tau}{(t - \tau)^{\mu}}~<~b
\Big(\int_0^{t}\frac{g(\tau) d \tau}{(t - \tau)^{\mu}}\Big)^2 $$
Lemma is proved.~$\blacktriangleleft$

 From  the basis equation (2,20) and the inequality (3,1), integrating over the domain $\Omega$, and, using
 the Holder inequality,  the estimates (2,9), (2,10) for
 Green function, we obtain the following estimate:
$$w(t)~~~<~~~f(t) ~~+~~~b~
\Big(\int_0^{t}\frac{w(\tau) + p(\tau)}{(t - \tau)^{\mu}}~d
\tau\Big)^2
 \eqno(3.2)$$where
$$w(\tau)~ =~\|\textbf{w}(x,\tau)\|_{L_2(\Omega)}~=~\sum_{i=1}^{i=3}\Big(\int_{\Omega}w_i^2(x,\tau)~dx\Big)^{1/2}
\geq 0;~$$
$$p(\tau) = \sum_1^3\Big\|\frac{\partial p(x,\tau}{\partial x_i}\Big\|_{L_2(\Omega)}~\geq 0;~~f(t) =
\|\textbf{f}(x,t)\|_{L_2(\Omega)} =
\sum_{i=1}^{i=3}\Big(\int_{\Omega}f_i^2(x,\tau)~dx\Big)^{1/2} \geq 0
\eqno(3.2')$$

\textbf{Remark 3-1.}~Let T > 0 be an arbitrary number. We have
proved that for all solutions $\textbf{w}(x,t), p(x,t)$ of the basis
equation  (2.20) the functions $w(\tau) =
\|\textbf{w}(x,\tau)\|_{L_2(\Omega)},~p(\tau)~ =~
\sum_1^3\Big\|\frac{\partial p(x,\tau}{\partial
x_i}\Big\|_{L_2(\Omega)}$ satisfy to the inequality (3.2). The
inequality (3.2) does not exclude the functions of the type~ $w(t) +
p(t) = \frac{t}{T-t} \nsubseteq L_2(0,T).$~ I.e. these functions
satisfy the inequality (3.2).~ Note that  $\mu: \frac{5}{8}\leq \mu
< \frac{3}{4}$ and $w(t) \in L_2(0,T).$ ~ By Proposition 3 and
estimate (2,17) in Proposition 8~ it follows that:
$$\Big\| \Big(\int_0^{t}\frac{(w(\tau) + p(\tau)) ~d \tau}{(t - \tau)^{\mu}}\Big)^2\Big\|_{L_2(0,T)} ~<~c\cdot
\|w(t) + p(t)\|_{L_2(0,T)} < c\cdot \|w(t)\|_{L_2(0,T)}
\eqno(3.3)$$In the Theorem 2.2 we assumed that for all possible
solutions $\textbf{w}(x,t)$ of the basic equation (2,20) the
function $w(t) = \|\textbf{w}(x,t) \|_{L_2(\Omega)} \in L_2(0,T).$

\textbf{By the basis equation (2,20)  we have proved the inequality
(3,2).~Below (see Theorem 3.1), we shall prove that for all
functions $w(t), p(t) \in L_2(0,T): \|p(t)\|_{L_2(0,T)} < c
\|w(t)\|_{L_2(0,T)} < \infty,$ ~satisfying the inequality (3,2), the
following a priori estimate: ~$\|w(t)\|_{L_2(0,T)} <
\sqrt{2}~\|f(t)\|_{L_2(0,T)}$~ holds.} $\blacktriangleleft$

\textbf{Lemma 3.2.}~From the estimates (2,17), (3,2) follows that:
$$w(t)~~~<~~~f(t) ~~+~~~b_1 \cdot
\int_0^{t}\frac{w^2(\tau)}{(t - \tau)^{\mu}}~d \tau
 \eqno(3.4)$$where $b_1 = b + c$ is constant independent on function
 w(t).

$\blacktriangleright$  By Holder inequality and the inequality
$(a+b)^2 < 2 a^2 + 2 b^2$ from the basic  inequality (3,2) follows
that:
$$w(t)~<~f(t) ~+~b
\Big(\int_0^{t}\frac{w(\tau) + p(t)}{(t - \tau)^{\mu/2}}\cdot
\frac{1}{(t-\tau)^{\mu/2}}~d \tau\Big)^2~<~f(t) +  \frac{2~ b }{1 -
\mu}\cdot T^{1-\mu} \cdot \int_0^{t}\frac{w^2(\tau)+ p^2(\tau)}{(t -
\tau)^{\mu}} d \tau $$
Using the estimate (2,17) in Proposition 8
and the inequality $(a_1 + a_2 + a_3)^2 < 3\cdot (a_1^2 + a_2^2 +
a_3^2)$, we infer:
$$ \int_0^t p^2(\tau)~ d \tau~=~ \int_0^t\Big(\sum_1^3\|p_{x_i}(x,\tau)\|_{L_2(\Omega)}\Big)^2 d \tau ~ <
~3 \cdot  \int_0^t\sum_1^3\|p_{x_i}(x,\tau)\|^2_{L_2(\Omega)}d \tau
<$$$$< 3 c \cdot \int_0^t \sum_1^3\|w_i(x,\tau)\|^2_{L_2(\Omega)} d
\tau < c
\int_0^t\Big(\sum_1^3\|w_{x_i}(x,\tau)\|_{L_2(\Omega)}\Big)^2 d
\tau~=~ 3 c \cdot \int_0^t w^2(\tau) d \tau $$ I.e.:~
$$w(t)~<~~f(t) ~+~b_1 \int_0^{t}\frac{w^2(\tau)+ p^2(\tau)}{(t -
\tau)^{\mu}} d \tau ;~~~~~~~~; \int_0^t p^2(\tau) d \tau~<~ c \cdot
\int_0^t w^2(\tau) d \tau    \eqno(3,4') $$ where ~$b_1~ =~ \frac{2~
b }{1 - \mu}\cdot T^{1-\mu}.$~~We shall prove by the second
inequality of (3,4') that there exists a constant c > 0:
$$\int_0^t \frac {p^2(\tau)}{(t - \tau)^{\mu}}~d \tau ~~<~~c\cdot \int_0^t \frac {w^2(\tau)}{(t - \tau)^{\mu}}~d
\tau    \eqno(3,4")$$We shall prove this estimate by contradiction
method and assumed that there exists a constants $c_n \rightarrow
\infty: $
$$\int_0^t \frac {p^2(\tau)}{(t - \tau)^{\mu}}~d \tau ~~>~~c_n\cdot \int_0^t \frac {w^2(\tau)}{(t - \tau)^{\mu}}~d
\tau$$Let us apply to this inequality the following
operator:~$J^{1-\mu}u(t) ~=~\int _0^t\frac {u(\tau)}{(t -
\tau)^{1-\mu}}~d \tau.$~~Then:~~$ \int_0^t p^2(\tau) d \tau~>~ c_n
\cdot \int_0^t w^2(\tau) d \tau .$~~But:~$\int_0^t p^2(\tau) d
\tau~<~ c \cdot \int_0^t w^2(\tau) d \tau.$~~And this contraction
proves the estimate (3,4").~Lemma 3,2 is proved.
$\blacktriangleleft$

\textbf{Remark 3-2.} Below, using the
 \textbf{Riccati's} replacement of the function w(t), from the estimate (3,4) we derive the
following estimate: $\|w(t)\|_{L_2(0,T)} <
\sqrt{2}~\|f(t)\|_{L_2(0,T)}$.~~For these aims we  consider the
following equation.

\textbf{The equation of Riccati.}~~In 1715, Riccati has studied the
following nonlinear equation on the segment [0,T] where T > 0 is an
arbitrary real number [10 p.41]:$$\frac{d z(t)}{d
t}~~=~~f(t)~+~b~z^2(t);~~~~z(0) = 0    $$ By the replacement of the
unknown function ~$z(t)=- \frac{1}{b}\cdot \frac{u'(t)}{u(t)},$~
this nonlinear equation is reduced to the following linear equation
of the second order:~~
$$\frac{d^2 u(t)}{d t^2}~+~b f(t) u(t)~=~0;~~~\frac{d u(t)}{d
t}|_{t=0} = 0.$$ Using the Riccati's result , we prove the following
key proposition. \vskip5mm \textbf{Theorem~3.1.}~For all functions
$w(t) \in L_2(0,T)$ ~satisfying  the inequality  (3,4) the following
estimate holds:
$$\|w(t)\|_{L_2(0,T)} < \sqrt{2}\cdot~\|f(t)\|_{L_2(0,T)}    \eqno(3,5)$$~~
This estimate does not depends on the number $b_1$ in (3,4).

$\blacktriangleright \blacktriangleright$~In (3.4) we  make the
replacement of the function:
$$w_1(t)~~~=~~~~ \int_0^{t}\frac{w^2(\tau))}{(t - \tau)^{\mu}}~d \tau      \eqno(3,6)$$
and, using the inequality $(a+b)^2 < 2 a^2 + 2 b^2,$ we rewrite the
basis inequality (3,4) as:
$$\int_0^t \frac{\frac{d w_1(\tau)}{d \tau}}{(t - \tau)^{1 - \mu}}
~d \tau~~<~~2 f^2(t)~+~ 2 b_1^2\cdot w_1^2(t) \eqno(3,7)$$Let~ $k :
0 < k < \infty$ and $s: 0 < s < 1$ - are an arbitraries real numbers
and for t > 1 we present the inequality (3,7) as:
$$\int_0^t \frac{\frac{d w_1(\tau)}{d \tau}}{(t - \tau)^{1 - \mu}}
~d \tau~<~2 f^2(t)~+~k \int_0^t \frac{w_1^2(\tau) d \tau}{(t -
\tau)^{1 - \mu}} ~+~2b_1^2\cdot t^{\mu}  w_1^2(t) ~-~k \int_0^t
\frac{w_1^2(\tau) d\tau}{(t - \tau)^{1 - \mu}} \eqno(3,8)$$ Applying
the Mellin transformation, for $s: 0 < s +\mu < 1$  we obtain:
$$\int_0^{\infty} t^{s-1}\cdot \Big(2 b_1^2\cdot t^{\mu}  w_1^2(t) - k \int_0^t \frac{w_1^2(\tau) d\tau}{(t
- \tau)^{1 - \mu}}\Big) ~d t=$$$$ = \Big(2 b_1^2 - k\cdot
\int_0^1\frac{d \tau}{\tau^{s+\mu} \cdot (1 -
\tau)^{1-\mu}}\Big)\cdot \int_0^{\infty} \tau^{s+\mu - 1}
w^2_1(\tau) d \tau$$~From this equality follows that: $\Big(2
b_1^2\cdot t^{1-\mu} w_1^2(t)  - k \int_0^t \frac{w_1^2(\tau)
d\tau}{(t - \tau)^{1 - \mu}}\Big) < 0$ for all numbers  $k \gg 1:$
$$ \Big(2 b_1^2 - k\cdot \int_0^1\frac{d \tau}{\tau^{s+\mu} \cdot (1 -
\tau)^{1-\mu}}\Big)~<~0 \eqno(3,9)$$ ~where~ $s: 0 < s +\mu < 1.$~
Using this inequality, we rewrite the inequality (3,8) as follows:
$$\int_0^t \frac{\frac{d w_1(\tau)}{d \tau}}{(t - \tau)^{1 - \mu}}
~d \tau~<~2 f^2(t)~+~k \int_0^t \frac{w_1^2(\tau) d \tau}{(t -
\tau)^{1 - \mu}}   \eqno(3,10)$$In this inequality , as Riccati , we
shall make the replacement of the function $w_1(\tau)$
$$w_1(\tau)~~=~~-~\frac{1}{k}~\frac{z'(\tau)}{z(\tau)};~~~~~;
\frac{d w_1(\tau)}{d \tau} ~=~- \frac{1}{
k}\cdot\frac{z''(\tau)}{z(\tau)} + \frac{1}{k}
\Big(\frac{z'(\tau)}{z(\tau}\Big)^2 \eqno(3,11)$$ From the
definition of the function $w_1(t)$ by (3,6) follows that:
$$z(t)~=~z(0)\cdot e^{- k \int_0^tw_1(\tau) d \tau}
~=~z(0)\cdot e^{- \frac{k}{1-\mu}\int_0^tw^2(\tau) (t -
\tau)^{1-\mu}d\tau} \eqno(3,11')$$ Since $w_1(0) = 0$, then $z'(0) =
\frac{d z(t}{d t}|_{t = 0} ~=~ 0.$~~From the inequality (3,10) we
have
$$ - \frac{1}{ k}~\int_0^t\frac{\frac{1}{z(\tau)}~\frac{d^2
z(\tau)}{d \tau^2}}{(t - \tau)^{1-\mu}}~d \tau ~~< ~2
f^2(t)\eqno(3,12)$$Or
$$ - \frac{1}{ k}~\int_0^t~  \frac{1}{z(\tau}\cdot \frac{d^2
z(\tau)}{d \tau^2}~d \tau ~~~<~~ \frac{2} {\Gamma_{1-\mu}^{\mu}}
\cdot \int_0^t \frac {f^2(\tau)}{(t - \tau)^{\mu}}~d \tau
\eqno(3,13)$$ Let us denote:
$$ \frac{1}{ k}~\int_0^t~  \frac{1}{z(\tau)}\cdot \frac{d^2
z(\tau)}{d \tau^2}~d \tau ~+~\frac{2}{\Gamma_{1-\mu}^{\mu}} \cdot
\int_0^t \frac {f^2(\tau)}{(t - \tau)^{\mu}}~d \tau~=~g(t)~>~0$$Then
$$ \frac{1}{ k}~  \cdot
\frac{d^2 z(t)}{d t^2}~~+~\frac{2}{\Gamma_{1-\mu}^{\mu}} \cdot z(t)
\frac{d}{d t} \int_0^t \frac {f^2(\tau)}{(t - \tau)^{\mu}}~d\tau
~=~z(t)~\frac{d g(t)}{d t} \eqno(3,13') $$Since $g(t) > 0 , g(0) =
0, -\frac{d z(\tau)}{d \tau} ~>~0$, integrating by parts, we obtain:
$$\int_0^t z(\tau)\cdot \frac{d}{d \tau} g(\tau~=~z(t)
g(t)~+~\int_0 ^t\Big(- \frac{d z(\tau)}{d \tau}\Big)\cdot g(\tau) d
\tau~~>~~0$$Since $\frac{d z(t)}{d t}|_{t = 0}~=~0$, integrating the
equation (3,13') over [0,t] , we obtain the following
$\textbf{important inequality}$:
$$-\frac{d z(t)}{d t} ~~~<~~k~\cdot \int_0^t z(\tau)~\frac{dF_{\mu}(\tau)}{d \tau} ~d \tau \eqno(3,14)$$ where
$$F_{\mu}(t)~~=~\frac{2}{\Gamma_{1-\mu}^{\mu}} \cdot
\int_0^t \frac {f^2(\tau)}{(t - \tau)^{\mu}}~d \tau
~~~~~\blacktriangleleft\eqno(3,15)$$

Integrating by parts,from the basis inequality (3,14), we have:
$$ - \frac{1}{k}\cdot \frac{d z(t)}{d t}~<~z(t)\cdot
F_{\mu}(t)~+~\int_0^t F_{\mu}(\tau) \Big(- \frac{d z(\tau)}{d
\tau}\Big)\cdot d \tau     \eqno(3,16)$$ In order to  proof  the a
priori estimate , it is \textbf{necessary} that the function
~$F_{\mu}(t) = \frac{2}{\Gamma_{1-\mu}^{\mu}}\cdot \int_0^t \frac
{f^2(\tau)}{(t - \tau)^{\mu}} d \tau$~ is increasing on [0,T].
~Otherwise,~the proof of this estimate  is difficult.

\textbf{Remark  3-3.}~Not for all positive  right-hand side $f^2(t):
f^2(t) \in L_1(0,T)$~ and real numbers  $\mu: 1/2 < \mu < 1$~ the
function  ~$f_{\mu}(t) = \frac{sin \pi \delta}{\pi} \int_0^t \frac
{f^2(\tau)}{(t - \tau)^{\mu}} d \tau$~ is increasing on [0,T].~For
example, we consider the following function ~$f_{\mu}(t) = \frac{sin
\pi \mu}{\pi} \int_0^t \frac {1-k\tau}{(t - \tau)^{\mu}} d
\tau$~where $(1-k \tau) > 0$ on (0,t), i.e. $0< \tau < 1/k.$~ Then
~$ \frac{d f_{\mu}(t)}{d t}~=~\frac{sin \pi \mu}{\pi} \cdot
\frac{1}{t^{\mu}}\cdot \Big(1 - \frac{k}{1-\mu}\cdot t \Big),$~ and
for ~$t:\frac{1-\mu}{k} <~t~< \frac{1}{k}$~ the function
$f_{\mu}(t)$ is decreasing on $(\frac{1-\mu}{k},
\frac{1}{k})$.$\blacktriangleleft$

Let $f(t): f(t) \in L_2(0,T)$ - is an arbitrary function. To define
the right-hand side $f(t),$~ for which the function $F_{\mu}(t)$ is
increasing on [0,T], we introduce the following functional spaces:
$$L_2^+ (0, T) ~~=~~\Big\{ f^2(t): f(t) > 0 ,~~\int_0^T f^2(\tau) d \tau ~<~ \infty,\Big\} $$
$$L_{1-\mu}^+ (0, T) ~~=~~\Big\{ f_{1-\mu}(t)~ : ~ f_{1-\mu}(t)~=~\int_0^t \frac {g^2(\tau)}{(t - \tau)^{1-\mu}}~d
\tau\Big\}  \eqno(3,17) $$where~$g(t):~g(t) \in L_2^+(0,T)$~is an
arbitrary function.~~

\textbf{Remark 3-4.} For all functions ~$\{f_{1-\mu}(t) = \int_0^t
\frac {g^2(\tau)}{(t - \tau)^{1-\mu}}~d \tau \in L_{1-\mu}^+(0,T)\}$
the functions $$F_{\mu}(t) = \frac{2}{\Gamma_{1-\mu}^{\mu}}\cdot
\int_0^t \frac{f_{1-\mu}(\tau)}{(t - \tau)^{\mu}} d \tau = \frac{2}
{\Gamma_{1-\mu}^{\mu}}\cdot \int_0^t \frac{d\tau}{(t - \tau)^{\mu}}
\int_0^{\tau}\frac{g^2(\tau_1) d \tau_1}{(\tau - \tau_1)^{1-\mu}} =
2 \int_0^t g^2(\tau) d \tau   \eqno(3,18)$$  are increasing on
[0,T]. This remark is $\textbf{important}$ for the proof the a
priori estimate. $\blacktriangleleft$

\textbf{Lemma 3.3}~For any right-side $f(t)
~=~f_{1-\mu}(t)~=~~\int_0^t \frac {g^2(\tau)}{(t - \tau)^{1-\mu}}~d
\tau \in L_{1-\mu}^+(0,T)$~in the basis equations (3,14) and (3,15)
the following a priori estimate holds:
$$\|w(t)\|_{L_2(0,T)} ~<~ \sqrt{2}\cdot ~\|f_{1-\mu}(t)\|_{L_2(0,T)}    \eqno(3,19)$$
$\blacktriangleright$~~ Since $\Big(-\frac{d z(t)}{d
~t}\Big)~>~0$~and for any function~$f_{1-\mu}(t)~=~ \in
L_{1-\mu}^+(0,T)$~ the function ~$F_{\mu}(t)~=~2\cdot \int_0^t
g^2(\tau) d \tau$~ is increasing (Remark 3.4), we rewrite the
inequality (3,16) as:
$$ - \frac{1}{k}\cdot \frac{d z(t)}{d t}~<~~z(t)\cdot
F_{\mu}(t)~+~F_{\mu}(t)\cdot \int_0^t  \Big(- \frac{d z(\tau)}{d
\tau}\Big)\cdot d \tau     \eqno(3,20)$$ Let us denote:
~$\int_0^t\Big(- \frac{d z(\tau)}{d \tau}\Big) d \tau~=~z_1(t).$~~~
Then this equation will accept the following kind: $$\frac{d
z_1(t)}{d t} - k\cdot F_{\mu}(t)z_1(t)~<~k \cdot F_{\mu}(t)\cdot
z(t)
 \eqno(3,20')$$~ Since $z_1(0) = 0$~and~$k\cdot F_{\mu}(t)\cdot z(t) > 0,$ from this inequality and
 Gronwall's Lemma  we infer that:
$$z_1(t) = \int_0^t\Big(- \frac{d z(\tau)}{d \tau}\Big) d \tau~<~
k\cdot \Big(\int_0^t F_{\mu}(\tau) z(\tau) e^{-k
\int_0^{\tau}F_{\mu}(\tau_1)d \tau_1} d \tau\Big)\cdot e^{k \int_0^t
F_{\mu}(\tau)d\tau}$$ Since $ z(t) > 0$ and the function $F_{\mu}(t)
> 0$ is  increasing, we rewrite this inequality:
$$z(0) - z(t)~<~ k \cdot F_{\mu}(t) \Big(\int_0^t z(\tau)
e^{-k \int_0^{\tau}F_{\mu}(\tau_1)d \tau_1} d \tau\Big)\cdot e^{k
\int_0^tF_{\mu}(\tau)d\tau} \eqno(3,21)$$ Let us denote:
$$z_2(t)~=~~\int_0^t z(\tau) e^{-k\cdot \int_0^{\tau}F_{\mu}(\tau_1)d
\tau_1} d \tau   \eqno(3,22)$$ ~and we present   (3,21) in the
following form:
$$\frac{d z_2(t)}{d t}~+~k \cdot
F_{\mu}(t) \cdot z_2(t)~~>~~z(0) \cdot  e^{-k
\int_0^{t}F_{\mu}(\tau)d \tau} \eqno(3,23)$$ Using the Gronwall's
Lemma and $z_2(0) = 0$ , we infer  that:
$$z_2(t)~~>~~z(0)\cdot t \cdot e^{- k  \int_0^{t}F_{\mu}(\tau)d
\tau}  \eqno(3,24)$$ Since~ $z_2(t) = \int_0^t z(\tau) e^{-k
\int_0^{\tau}F_{\mu}(\tau_1)d \tau_1} d \tau~<~\int_0^t z(\tau)d
\tau,$  it follows from (3,22),(3,24) that
$$z(0)\cdot t~~~<~~ \Big(\int_0^t z(\tau)  d \tau \Big)\cdot e^{ k
\int_0^{t}F_{\mu}(\tau)d \tau} \eqno(3,25)$$ where the function z(t)
is defined by (3,11') and (3,6):~
$$z(t)~=~z(0)\cdot e^{- k \int_0^tw_1(\tau) d \tau}
~=~z(0)\cdot e^{- \frac{k}{1-\mu}\int_0^tw^2(\tau) (t -
\tau)^{1-\mu}d\tau} \eqno(3,25')$$

\textbf{Remark 3.5}.~From (3,25') we have:
$$\frac{d z(t)}{dt} = - z(0)\cdot k\cdot \Big(\int_0^t
\frac{w^2(\tau)}{(t-\tau)^{\mu}}d \tau\Big) \cdot e^{-
\frac{k}{\mu}\int_0^tw^2(\tau) (t - \tau)^{1-\mu}d\tau} < 0.$$ If
for a some~$t_0 \in (0,\infty)$~ $\frac{d z(t}{d t}\big|_{t = t_0} =
0$ , then the following two cases are possible:
\textbf{1.}$\int_0^{t_0} \frac{w^2(\tau)}{(t_0-\tau)^{\mu}}d \tau
~=~0$,~~or~~\textbf{2.}$\int_0^{t_0}w^2(\tau) (t_0 -
\tau)^{1-\mu}d\tau~=~\infty.$~~If \textbf{1.}$\int_0^{t_0}
\frac{w^2(\tau)}{(t_0-\tau)^{\mu}}d \tau ~=~0,$~then~ $w(t) \equiv
0$ ~on $[0,t_0),$~since $w(t) \geq 0$ on $[0, t_0)$.~~If~
\textbf{2.}$\int_0^{t_0}w^2(\tau) (t_0 -
\tau)^{1-\mu}d\tau~=~\infty,$~then below we shall prove that for all
positive functions $w^2(t)$ satisfying to the basis inequality
(3,25): ~~$\int_0^{t_0}w^2(\tau) (t_0 -
\tau)^{1-\mu}d\tau~<~\infty.$ $\blacktriangleleft$

\textbf{The following is the key Lemma for the proof of Theorem
2.2.}

\textbf{Lemma 3,5'.} Let t > 0 is an arbitrary real number.~For all
positive functions ~$w^2(\tau)$ satisfying to the basis inequality
(3,25) follows that:~$$\int_0^{t}w^2(\tau) (t -
\tau)^{1-\mu}d\tau~<~\infty        \eqno(3,36)$$ where~$\mu: 1/2 <
\mu < 1.$~~Or, passing to limit ~$\mu \rightarrow 1$,~we obtain:
$$\int_0^t w^2(\tau)~d \tau~~~<~~\infty          \eqno(3,36')$$

$\blacktriangleright$~We shall prove this Lemma by contradiction
method and  rewrite the basis inequality (3,25):
$$z(0)\cdot t~~~<~~ \Big(\int_0^t z(\tau)  d \tau \Big)\cdot e^{ k
\int_0^{t}F_{\mu}(\tau)d \tau} \eqno(3,25)$$ where the function z(t)
is defined by (3,11'):~
$$z(t)~=~z(0)\cdot e^{- \frac{k}{1-\mu}\int_0^tw^2(\tau) (t -
\tau)^{1-\mu}d\tau} \eqno(3,25')$$ Let for a some number $t_0:$
~$\int _0^{t_0} w^2(\tau)(t_0 - \tau)^{1-\mu}~d \tau
~=~\infty.$~I.e.~$w^2(\tau)~\approx \frac{c}{|t_0 - \tau|^{\lambda
}\cdot |t_0 - \tau|^{1-\mu}}$ ~~where ~$\lambda \geq 1$ is a real
number.~Then for all $t \geq t_0$:~~$\int _0^{t} w^2(\tau)(t -
\tau)^{1-\mu}~d \tau ~=~\infty.$ ~~ And ~~~$z(t)~=~z(0)\cdot e^{-
\frac{k}{1-\mu}\int_0^tw^2(\tau) (t - \tau)^{1-\mu}d\tau}~\equiv
0$~for~$t \geq t_0$.~~Since $e^{- \frac{k}{1-\mu}\int_0^tw^2(\tau)
(t - \tau)^{1-\mu} \tau}\cdot t\Big|_0^{t_0}~=~0$, then, integrating
by part, we have:
$$\frac{1}{z(0)}\int_0^{t_0} z(\tau)  d \tau~=~+~\int_0^{t_0}t\cdot
e^{-\frac{k}{1-\mu}\int_0^tw^2(\tau) (t - \tau)^{1-\mu}d\tau}~d
\Big( \frac{k}{1-\mu}\int_0^tw^2(\tau) (t - \tau)^{1-\mu}
\tau\Big)$$In this equality we shall make the replacement of
variable $\frac{k}{1-\mu}\int_0^tw^2(\tau) (t - \tau)^{1-\mu} \tau =
t_1.$ Then:
$$\frac{1}{z(0)}\int_0^{t_0} z(\tau)  d \tau~<~t_0\cdot
\int_0^{\infty}e^{-t_1}d~ t_1~~=~~t_0$$ From the definition of the
function $F_{\mu}(t)$ ~by (3,15) we obtain:~
$$\int_0^t F_{\mu}(\tau) d \tau~~=~\frac{2}{\Gamma_{1-\mu}^{\mu}} \cdot
\frac{1}{1-\mu}\cdot \int_0^t f^2(\tau)(t - \tau)^{1-\mu}~d \tau
\eqno(3,15')$$Using the following limits: ~$lim_{\mu \rightarrow
1}\Gamma(1-\mu)\cdot (1-\mu)~=~1$~and~$lim_{\mu\rightarrow 1}(t -
\tau)^{1-\mu} ~=~1$, we get: ~$\lim _{\mu\rightarrow 1}\int_0^t
F_{\mu}(\tau) d \tau~=~\int_0^t f^2(\tau) d \tau$.~Using these
facts, from the basis inequality (2,25) we have: ~~$$ t~~<~~ t_0
~~\cdot~~ e^{ \int_0^{T}f^2(\tau)d \tau}$$ But, for $t \gg 1$ ~we
have received the contradiction.~ Therefore~$\int_0^{t}w^2(\tau) (t
- \tau)^{1-\mu}d\tau~<~\infty. $~

Lemma is proved.$\blacktriangleleft$

\textbf{Remark 3.35'.}~Hence,  for \textbf{any positive} function
$w^2(t) \geq 0$~satisfying to the basis inequality (3,25) the
function ~$z(t) > 0$ ~is continuous on $[0, \infty)$ and
monotonously decreasing from z(0) to zero on
$[0,\infty)$.~$\blacktriangleleft$

Let T > 0 be an arbitrary real number. We rewrite the basis
inequality (3,25) for t = T  as:
$$T\cdot z(0)~~~<~~ \Big(\int_0^{ T} z(t) d t~+~\int_T^{2 T} z(t) d t \Big)\cdot e^{ k
\int_0^T F_{\mu}(\tau)d \tau} \eqno(3,27)$$ Since the function
~$z(t)~=~z(0)\cdot e^{-\frac{k}{1-\mu}\int_0^tw^2(\tau) (t -
\tau)^{1-\mu}d\tau}$~is continuous on $[0, \infty)$ and monotonously
decreasing from z(0) to zero, there exists the numbers~$t_1 , t_2$:
$$t_1~:~0~<~t_1~<~T~;~~~~~~~~;~t_2:~T~<~t_2~<~2 T    \eqno(3,27')$$
such that the following equalities holds:
$$\int_0^{T}  e^{-  \frac{k}{1-\mu} \int_0^tw^2(\tau) (t -
\tau)^{1-\mu}d\tau}  d \tau~=~ \Big(e^{- \frac{k}{1-\mu}
\int_0^{t_1}w^2(\tau) (t - \tau)^{1-\mu}d\tau}\Big)\cdot T
\eqno(3,28)$$
$$\int_T^{2 T}  e^{- \frac{k}{1-\mu} \int_0^tw^2(\tau) (t -
\tau)^{1-\mu}d\tau}  d \tau~=~ \Big(e^{- \frac{k}{1-\mu}
\int_0^{t_2} w^2(\tau) (t - \tau)^{1-\mu}d\tau}\Big)\cdot T
\eqno(3,28')$$ Using these equalities, we rewrite the inequality
(3,27) as:
$$1~<~~\Big(e^{-  w_2(t_1)}~~+~~e^{-  w_2(t_2)}\Big)\cdot
~e^{k\cdot \int_0^{T} F_{\mu}(\tau)~d \tau} \eqno(3,29)$$ where
$$w_2(t_1)~~=~~
\frac{k}{1-\mu} \int_0^{t_1}w^2(\tau) (t - \tau)^{1-\mu} d
\tau;~~~~~;w_2(t_2)~~=~~ \frac{k}{1-\mu} \int_0^{t_2}w^2(\tau) (t -
\tau)^{1-\mu} d \tau \eqno(3,29')$$and present the inequality (3,29)
as:$$ e^{w_2(t_2)}~~<~~ \Big(\frac{e^{ w_2(t_2)}}{e^{w_2(t_1)}}
~+~1\Big)\cdot e^{k\cdot \int_0^{T}F_{\mu}(\tau) d \tau}
\eqno(3,30)$$ Let in formula (3,15)~the right-side
~$$f^2(t)~=~f_{1-\mu}^2(t)~=~\int_0^t \frac{g^2(\tau)}{(t -
\tau)^{1-\mu}} d \tau ~ ~\in~ L_{1-\mu}^+(0,T)
   \eqno(3,31)$$where the~ $g(t): ~g (t) ~\in~L_2(0,T)$~is an arbitrary
function.~~~Then,~from (3,15) we have:
$$F_{\mu}(t) ~ ~=~\frac{2}{\Gamma_{1-\mu}^{\mu}}\cdot
\int_0^t \frac {f^2_{1-\mu}(\tau)} {(t - \tau)^{\mu}} d \tau ~=
\eqno(3,32)$$$$=~2\cdot \frac{1}{\Gamma_{1-\mu}^{\mu}}\cdot \int_0^t
\frac {1} {(t - \tau)^{\mu}} \int_0^{\tau} \frac{g^2(\tau_1)~ d
\tau_1}{(\tau - \tau_1)^{1-\mu}}~d \tau~=~2\cdot \int_0^t
g^2(\tau)~d \tau$$From (3,31) and (3,32) follows that:
$$ g^2(t)~=~\frac{d}{d t} \int_0^t \frac{f_{1-\mu}^2(\tau)}{(t -
\tau)^{\mu}} d \tau;~~~~~~; F_{\mu} (t)~=~2 \int_0^tg^2(\tau) d
\tau~=~ 2\cdot \int_0^t \frac{f_{1-\mu}^2(\tau)}{(t - \tau)^{\mu}}~d
\tau \eqno(3,32')$$and from (3,32') we have:
$$\int_0^t F_{\mu}(\tau)~ d \tau ~~=~~\frac{2}{1-\mu}\cdot  \int_0^t (t -
\tau)^{1 - \mu}\cdot f_{1-\mu}^2(\tau)~~d \tau     \eqno(3,33)$$
Since the function~ $e^{- w_2(t)}$~is decreasing on $(0,\infty)$
~and~$t_2~>~t_1$, we have: ~$e^{ w_2(t_2)}/e^{w_2(t_1)}~<~1.$ And
from (3,30) follows that
$$e^{\frac{k}{1-\mu}\int_0^{t_2}w^2(\tau) (t - \tau)^{1-\mu}d\tau}~~<~~2\cdot
e^{2\cdot \frac{k}{1-\mu}\int_0^{T}f_{1-\mu}^2(\tau) (t -
\tau)^{1-\mu}d\tau}$$ or
$$k \int_0^{t_2}w^2(\tau) (t - \tau)^{1-\mu}d\tau~~<~~(1 - \mu)\cdot ln2~~+~
2 k \cdot \int_0^{T}f_{1-\mu}^2(\tau) (t - \tau)^{1-\mu}d\tau
\eqno(3,34)$$ Passing to the limit ~$\mu ~\rightarrow~1$~,~ from
~$t_2~>~T,$~$\lim_{\mu\rightarrow 1}(t-\tau)^{1-\mu}=1$ and this
inequality we obtain:
$$\int_0^T w^2(\tau) d\tau ~~~<~~2\cdot \int_0^T f^2_{1-\mu}(\tau) d\tau
\eqno(3,34')$$ Lemma 3.3 ~is \textbf{proved.} $\blacktriangleleft$

\textbf{The proof of Theorem 3.1.} ~~~ \textbf{Lemma 3-4.}~The space
of functions  $L_{1-\mu}^+(0,T)$ ~is dense in the space of functions
$L_2^+(0,T)$~ in the norm of the space $L_2(0,T)$.~~I.e. for all
functions  $f(t) \in L_2^+(0,T),$~ there exists a sequence of
functions:~~ $f_{n_{1-\mu}}(t)~\in
L_{1-\mu}^+(0,T):$~~$$lim_{n\rightarrow \infty}~\Big(\|
f_{n_{1-\mu}}(t)\|_{L_2(0,T)} - \|
f(t)\|_{L_2(0,T)}\Big)~<~lim_{n\rightarrow \infty} \|
f_{n_{1-\mu}}(t) ~-~f(t)\|_{L_2(0,T)}~~\rightarrow~~0 \eqno(3,35)$$
or
$$lim_{n\rightarrow \infty}~\Big(\|
f_{n_{1-\mu}}^2(t)\|_{L_1(0,T)} - \|
f^2(t)\|_{L_1(0,T)}\Big)~<~lim_{n\rightarrow \infty} \|
f_{n_{1-\mu}}^2(t) ~-~f^2(t)\|_{L_1(0,T)} \rightarrow 0
\eqno(3,35')$$ \textbf{Remark 3.6.} Let us note that for any n =
1,2,...the functions $f_{n_{1-\mu}}^2(t)~\in~L_{1-\mu}^+(0,T).$ From
Remark 3.4 follows that the functions $F_{n_{\mu}}(t) =
\frac{2}{\Gamma_{1-\mu}^{\mu}}\cdot \int_0^t
\frac{f_{n_{1-\mu}}(\tau)}{(t - \tau)^{\mu}} d \tau$~are increasing
on [0,T].~~And from (3,34) in Lemma 3.3 follow the following
estimates: $\|w(t)\|_{L_2(0,T)} ~<~ 2\cdot
~\|f_{n_{1-\mu}}(t)\|_{L_2(0,T)}.$ $\blacktriangleleft$

$\blacktriangleright$ We prove Lemma 3.4 by the contradiction
method.~ Let there exists a function $f_0(t):~f_0(t)\geq 0, f_0(t)
\neq 0,~ f_0(t) \in L_2(0,T),$ such that for all functions ~$f(t):
f(t)\geq 0,~~f(t) \in L_2(0,T)$~ the following equality is valid
$$\int_0^T f_0(t) \cdot \int_0^t \frac{f(\tau) d
\tau}{(t-\tau)^{1-\mu}} d t~~=~~\int_0^T f(\tau)\cdot
\Big(\int_{\tau}^T \frac{f_0(t) d t}{(t - \tau)^{1-\mu}}\Big)~d
\tau~=~0$$Since  ~$f(\tau) \geq 0$  is an arbitrary function, by
this equality it follows that for all $\tau \in [0,T]:$
~$\int_{\tau}^T \frac{f_0(t) d t}{(t -
\tau)^{1-\mu}}~\equiv~0.$~Then $f_0(t) \equiv 0$.~~~Lemma is proved.
$\blacktriangleleft$

Let $\epsilon_n : 0 < \epsilon_n \ll 1,~\lim _{n \rightarrow
\infty}\epsilon_n = 0$~are an arbitraries  real numbers.~Then for
any function $f(t) \in L_2(0,T),~f(t) > 0$ by (3,35') and for ~$n
\gg 1$ it follows that:
$$f^2(t) = f^2(t) - f_{n_{1-\mu}}^2(t) + f_{n_{1-\mu}}^2(t) <
\mid f^2(t) - f_{n_{1-\mu}}^2(t)\mid + f_{n_{1-\mu}}^2(t)~ <~
\epsilon_n ~+~f_{n_{1-\mu}}^2(t)$$ Then:
$$F_{\mu}(t)~=~\frac{2}{\Gamma_{1-\mu}^{\mu}}\cdot \int_0^t
\frac{f^2(\tau)}{(t
-\tau)^{\mu}}d\tau~<~~\frac{2}{\Gamma_{1-\mu}^{\mu}} \cdot \int_0^t
\frac { \epsilon_n + f_{n_{1-\mu}}^2(\tau)}{(t - \tau)^{\mu}}~d \tau
~=~$$$$ ~=~\frac{2}{\Gamma_{1-\mu}^{\mu}}\cdot \Big(
\frac{t^{1-\mu}}{1-\mu} \cdot \epsilon_n~+~ \int_0^t \frac
{f_{n_{1-\mu}}^2(\tau}{(t - \tau)^{\mu}}~d \tau\Big)$$ Let us
denote:~
$$F_{\mu}^n(t)~~=~~\frac{2}{\Gamma_{1-\mu}^{\mu}}\cdot~\Big(\frac{t^{1-\mu}}{1-\mu} \cdot
\epsilon_n~+~ \int_0^t \frac {f_{n_{1-\mu}}^2(\tau) }{(t -
\tau)^{\mu}}~d \tau\Big) \eqno (3,36)$$ Since the functions
$f_{n_{1-\mu}}^2 \in L_{1-\mu}^+(0,T)$ ~and the function
~$t^{1-\mu}$~is increasing, the functions $F_{\mu}^n(t)$ are
increasing on [0,T] ~and~$F_{\mu}(t)~<~F_{\mu}^n(t).$~~Therefore,
from the basis inequality (3,25) we have:
$$z(0)\cdot t~~~<~~ \Big(\int_0^t z(\tau)  d \tau \Big)\cdot e^{ k
\int_0^{t}F_{\mu}^n(\tau)d \tau} \eqno(3,25')$$ and
$$\int_0^t F_{\mu}^n(\tau) d
\tau~=~\frac{2}{\Gamma_{1-\mu}^{\mu}}\cdot
\Big(\frac{t^{2-\mu}}{(1-\mu)(2-\mu)}\epsilon_n~+~\frac{1}{1-\mu}\cdot
\int_0^t(t-\tau)^{1-\mu}    f_{n_{1-\mu}}^2(\tau) d \tau \Big)$$
Further similarly to the proof of an estimate (3,34), we obtain:
$$k \int_0^{t_2}w^2(\tau) (t - \tau)^{1-\mu}d\tau~~<~~(1 - \mu)\cdot ln2~~+~
 \int_0^T F_{\mu}^n(\tau) d \tau \eqno(3,34')$$
Using the following limits: ~$lim_{\mu \rightarrow 1} \Gamma(1 -
\mu)\cdot (1-\mu)~=~1$,~~$lim_{n\rightarrow \infty}~\|
f_{n_{1-\mu}}^2(t)\|_{L_1(0,T)} = \|
f^2(t)\|_{L_1(0,T)},$~$lim_{n\rightarrow \infty} \epsilon_n ~=~0$
and Remark 3.6 , ~passing to the limits~ $\mu \rightarrow 1$~and~$n
\rightarrow \infty,$~ from the inequality  (3,34') we get:
$$\int_0^T w^2(\tau) d\tau ~~~<~~2\cdot \int_0^T f^2(\tau) d\tau
\eqno(3,37)$$~~Theorem 3,1 is proved.~$\blacktriangleleft
\blacktriangleleft$

\begin{center}\textbf{Proof of Theorem 2,1.} \end{center}
\textbf{Definition 3.} If the sequence of vector-functions
~$\{\textbf{w}_n(x,t)\}$~ $\textbf{weakly}$ converges to the
vector-function ~$\textbf{w}_0(x,t)$ in the space
$\textbf{L}_2(Q_t),$~ then we denote : ~$\textbf{w}_n(x,t)
\rightharpoonup \textbf{w}_0(x,t).$~I.e. for an arbitrary
vector-function ~$\textbf{u}(x,t) \in \textbf{L}_2(Q_t)$~ the
following convergence is valid:~ $\big(
\textbf{w}_{n}(x,t)~,~\textbf{u}(x,t)\big)_{L_2(Q_t)}~\rightarrow~\big(
\textbf{w}_0(x,t)~,~\textbf{u}(x,t)\big)_{L_2(Q_t)}~as~n \rightarrow
\infty$

$\textbf{2)}$~ If the sequence of vector-functions
~$\{\textbf{w}_n(x,t)\}$~ $\textbf{strongly }$ converges to the
vector-function ~$\textbf{w}_0(x,t)$ in the space
$\textbf{L}_2(Q_t)$ ,~ then we denote: ~$\textbf{w}_n(x,t)
\Rightarrow \textbf{w}_0(x,t).$~I.e. $lim_{n\rightarrow
\infty}\|\textbf{w}_n(x,t) - \textbf{w}_0(x,t)\|_{L_2(Q_t)}
\rightarrow 0.$

\textbf{Lemma 3.4}~Let the vector-function~~$\textbf{w}(x,t)~=~(w^i
(x,t))_{i=1,2,3} \in~\textbf{L}_2(Q_t)$.~~~ I.e.
$\|\textbf{w}(x,\tau)\|_{L_2(Q_t)}~=~\sum_{i=1}^{i=3}(\int_{Q_t}(w^i)^2(x,\tau)~dx
d \tau)^{1/2}~<~\infty.$~We define the following nonlinear operator
K on the vector-space $\textbf{L}_2(Q_t)$:
$$ K\ast \Big(\textbf{w}(x,t)\Big)~=~G w^j(x,t) G_{x_j}\textbf{w}(x,t) = \Big(\sum_{j=1}^3
G w^j(x,t) G_{x_j} w^i(x,t)\Big)_{i=1,2,3}$$ Let us prove that  ~
$$K\ast \textbf{w}_{n}(x,t) \Rightarrow K\ast
\textbf{w}_0 (x,t)~   \eqno(3.38)$$ as ~ $\textbf{w}_n(x,t)
\rightharpoonup \textbf{w}_0 (x,t)$.~

It follows from this proposition  that the operator K is compact on
the vector-space $\textbf{L}_2(Q_t)$.~It is follows from book [3
p.42]

$\blacktriangleright$~Then
$$\|K\ast \textbf{w}_{n}(x,t) - K\ast \textbf{w}_0(x,t)\|_{L_2(Q_t)}
\leq c \sum_{i,j =1}^3 \|G w_{n}^j(x,t) G_{x_j}w_{n}^i(x,t) - G
w_{0}^{j}(x,t) G_{x_j}w_{0}^{i}(x,t)\|_{L_2(Q_t)}$$ Let us estimate
a each member: $$\Big\|G w_{n}(x,t) G_{x}w_{n}(x,t) - G w_{0}(x,t)
G_{x}w_{0}(x,t)\Big\|_{L_2(Q_t)}\leq    \eqno(3.39)$$
$$\leq \Big\|G_x w_{n}(x,t)\Big(G w_{n}(x,t) - G w_{0}(x,t)\Big)
\Big\|_{L_2(Q_t)} + \Big\|G w_{0}(x,t) \Big(G_{x}w_{n}(x,t) - G
_{x}w_{0}(x,t)\Big)\Big\|_{L_2(Q_t)}$$ We
 obtain the following  estimates, using the formula  (2.12.2)
in Proposition 5, the inequality (2,4) in Proposition 1 and
Proposition 3 for ~$\mu: 5/8 < \mu < 1,$:
$$\Big\|G_x w_{n}(x,t)\Big(G w_{n}(x,t) - G w_{0}(x,t)\Big)
\Big\|_{L_2(Q_t)}  \leq$$$$\leq c~\Big\|G_x
w_n(x,t)\Big\|_{L_{\frac{6}{5-4\mu}, \frac{2}{2\mu - 1}}(Q_t)}\cdot
\Big\|G \Big(w_{n}(x,t) - w_{0}(x,t)\Big)\Big\|_{L_{\frac{3}{2 \mu -
1},\frac{1}{1 - \mu}}(Q_t)} \leq$$$$\leq
c~\Big\|w_n(x,t)\Big\|_{L_2(Q_t)} \cdot \Big\|G \Big(w_{n}(x,t) -
w_{0}(x,t)\Big)\Big\|_{L_{\frac{3}{2 \mu - 1},\frac{1}{1 -
\mu}}(Q_t)}~~\rightarrow~0   \eqno(3,40)$$ We obtain the following
estimates , using the formula (2.12.1) in Proposition 5 and
  the inequality (2,4) in Proposition 1 for ~$\mu: 5/8 < \mu < 1,$:
$$\Big\|G w_{0}(x,t)\Big(G_x w_{n}(x,t) - G_x w_{0}(x,t)\Big)
\Big\|_{L_2(Q_t)}  \leq$$$$\leq c~\Big\|G
w_0(x,t)\Big\|_{L_{\frac{6}{3-4\mu}, \frac{2}{2\mu - 1}}(Q_t)}\cdot
\Big\|G_x \Big(w_{n}(x,t) - w_{0}(x,t)\Big)\Big\|_{L_{\frac{3}{2 \mu
},\frac{1}{1 - \mu}}(Q_t)} \leq$$$$\leq
c~\Big\|w_0(x,t)\Big\|_{L_2(Q_t)} \cdot \Big\|G_x\Big(w_{n}(x,t) -
w_{0}(x,t)\Big)\Big\|_{L_{\frac{3}{2 \mu},\frac{1}{1 -
\mu}}(Q_t)}~~\rightarrow~0   \eqno(3,41) $$ Note that by the weakly
convergence  $w_n(x,t)-w_0(x,t) \rightharpoonup 0$~ it follows that
for any number n there is a constant c such that:
$\|w_n(x,t)\|_{L_2(Q_t)} < c \|w_0(x,t)\|_{L_2(Q_t)}.$~ Since $G
\Big(w_{n}(x,t) - w_{0}(x,t)\Big) \in W_2^{2,1}(Q_t)$,  the space
$W_2^{2,1}(Q_t)$~ is compactly enclosed into the space
$L_{\frac{3}{2\mu-1},\frac{1}{1-\mu}}(Q_t)$~[5 p.78], and on the
space $L_{\frac{3}{2\mu-1},\frac{1}{1-\mu}}(Q_t)$~the compact
operator G translates the weekly convergence $w_n(x,t)-w_0(x,t)
\rightharpoonup 0$ to  the strongly convergence, then the strongly
convergence  (3.40)~is valid.

Since  $G_x \Big(w_{n}(x,t) - w_{0}(x,t)\Big) \in W_2^{1}(Q_t)$, the
space  $W_2^{1}(Q_t)$~is compactly enclosed into the space
$L_{\frac{3}{2\mu},\frac{1}{1-\mu}}(Q_t)$~[5 p.78], and on the space
$L_{\frac{3}{2\mu},\frac{1}{1-\mu}}(Q_t)$~the compact operator $G_x$
translates the weakly convergence $w_n(x,t)- w_0(x,t)
\rightharpoonup 0$ to  the strongly convergence, then the strongly
convergence (3.41) is valid . The strongly convergence (3,38)
follows by (3,39), (3,40), (3,41). Lemma is
proved.$\blacktriangleleft$

\textbf{The proof of Theorem 2.1.}~Let $\textbf{z}(x,t)
~=~(z_i(x,t))_{i=1,2,3} \in \textbf{L}_2(Q_t)$ be an arbitrary
vector-function.~~And the sequence vector-functions
$\textbf{w}^n(x,t)$ weakly converges to the vector-function
 $\textbf{w}^0(x,t),$~i.e.~$\textbf{w}^n(x,t) \rightharpoonup
\textbf{w}^0(x,t)$. On the vector-space $\textbf{L}_2(Q_t)$ we
define the following nonlinear operator $K_p$:
$$ K_p\ast \Big(\textbf{w}(x,t)\Big)~=~
\Big(\sum_{j=1}^{3}G \Big(w_j(x,t) + \frac{\partial p(x,t)}{\partial
x_j}\Big)\cdot G_{x_j}\Big(w_i(x,t) + \frac{\partial
p(x,t)}{\partial x_i}\Big)\Big)_{i=1,2,3}  \eqno(3,42)$$where the
functions $\frac{d p(x,t)}{d x_i}$ are defined by the functions
$w_i(x,t)$ from the formula  (2,16) in Proposition 8.~~ On the
vector space $\textbf{w}(x,t) ~=~\Big(w_i(x,t)\Big)_{i = 1,2,3}
~\in~\textbf{L}_2(Q_t)$ by the formula (2,18) in Proposition  8 we
have defined the following linear and bounded  operator:
$$P \ast \Big(w_i(x,t)\Big)_{i = 1,2,3}~=~\Big(\frac{d p(x,t)}{d
x_i}\Big)_{i = 1,2,3}~\in ~\textbf{L}_2(x,t) \eqno(3,43)$$ where the
functions  $\frac{d p(x,t)}{d x_i}$~ is defined by the functions
$w_i(x,t)$ ~from the formula (2,16) in Proposition 8. ~Since P is
linear and bounded operator on $\textbf{L}_2(Q_t)$, there exists the
linear and bounded connected  operator $P^{\ast}.$~ Let
$\textbf{z}(x,t) \in \textbf{L}_2(Q_t)$ is an arbitrary
vector-function and $\textbf{w}^n(x,t) \rightharpoonup
\textbf{w}^0(x,t)$.~Then:
$$\Big(\Big(w^n_i(x,t) + \frac{d p(x,t)}{d x_i}\Big)_{i=1,2,3},
\textbf{z}(x,t)\Big)_{\textbf{L}_2(Q_t)} = \Big(\big(w^n_i(x,t) +
P\ast(w_i^n(x,t)\big)_{i=1,2,3}~ ,~
\textbf{z}(x,t)\Big)_{\textbf{L}_2(Q_t)} =$$
$$= \Big(\textbf{w}^n(x,t)~ ,~ \textbf{z}(x,t) + P^{\ast}\ast
\textbf{z}(x,t)\Big)_{\textbf{L}_2(Q_t)} \rightarrow
\Big(\textbf{w}^0(x,t)~ ,~ \textbf{z}(x,t) + P^{\ast}\ast
\textbf{z}(x,t)\Big)_{\textbf{L}_2(Q_t)} = $$
$$= \Big(\textbf{w}^0(x,t) +
P\ast \textbf{w}^0(x,t)~ ,~
\textbf{z}(x,t)\Big)_{\textbf{L}_2(Q_t)}$$as $n \rightarrow \infty$.
~~I.e, it is proved that:  $\textbf{w}^n(x,t) + P\ast
\textbf{w}^n(x,t) ~\rightharpoonup \textbf{w}^0(x,t) + P\ast
\textbf{w}^0(x,t).$ Then, similarly to the proof in Lemma 3.4,
proves that:~ $K_p\ast \textbf{w}^{n}(x,t) \Rightarrow K_p\ast
\textbf{w}^0 (x,t).$ ~~Hence, it follows that on the vector-space
$\textbf{L}_2(Q_t)$ the nonlinear operator $K_p$ is compact.~[3
p.42]. ~~Therefore, it follows by the Leray -Schauder's theorem in
Proposition 6, that the basis equation (2,20) has at least one
solution  $\textbf{w}(x,t) \in \textbf{L}_2(Q_t)$ and it follows
from Theorem 3.1 that:
 $\|\textbf{w}(x,t)\|_{L_2(Q_T)} \leq \sqrt{2}\cdot\|\textbf{f}(x,t)\|_{L_2(Q_T)}.$
 ~~Then, it follows by Proposition 8 that there exists the smooth solution
  $\textbf{u}(x,t)~\in \textbf{W}_2^{2,1}(Q_T)\bigcap \textbf{H}_2(Q_t).$
 But the Navier-Stokes problem has the unique smooth solution [3p.139]. Therefore,
 $\textbf{w}(x,t) \in \textbf{L}_2(Q_T)$ is the unique solution to (2,20).
 The existence and smoothness of the solution to Navier-stokes equation is proved.$\blacktriangleleft$

\vskip 2mm

\textbf{4.~The Navier-Stokes problem for the inhomogeneous boundary
condition.} \vskip 2mm

Let~$\Omega \subset R^3$~be a finite domain bounded by a Lipschitz
surface ~$\eth \Omega$~and ~$Q_{T} = \Omega \times[0,T] , S = \eth
\Omega \times [0,T]~,~x =(x_1,x_2,x_3)$~and~$\textbf{u}(x,t)~=~
(u_i(x,t)_{i = 1,2,3}~,$ $\textbf{f}(x,t) = (f_i(x,t)_{i =1,2,3}$
are  vector-functions.~ Here~T > 0 is an arbitrary real number.~ The
Navier-Stokes equations are given by:
$$\frac{\partial u_i(x,t)}{\partial t}
~-~\rho~\triangle~u_i(x,t)~-~\sum_{j=1}^{3}u_j(x,t)~ \frac{\partial
u_i(x,t)}{\partial x_j}~+~\frac{\partial p(x,t)}{\partial
x_i}~~=~f_i(x,t)~\eqno(4.1),$$
$$div~\textbf{u}(x,t)~=~\sum_{i=1}^3\frac{\partial u_i(x,t)}{\partial x_i}~=~
0~~,i~=~1,2,3$$\textbf{The Navier-Stokes problem~1.} Find a
vector-function  $\textbf{u}(x,t)~=~(u_i(x,t))_{i=1,2,3} : \Omega
\times [0,T] \rightarrow R^3,$~the   scalar function $p(x,t):
 \Omega \times [0,T] \rightarrow R^1$  satisfying  the equation
(4.1) and   the following initial condition
$$\textbf{u} (x,0) =~\textbf{a}(x)~~,~~~~
 ~ \textbf{u}(x,t)\mid_{\partial \Omega \times [0,T]}~=~0
\eqno(4.2)$$where~~div $\textbf{a}(x) = \frac{d a_1(x)}{d x_1} +
\frac{d a_2(x)}{d x_2} + \frac{d a_3(x)}{d x_3} =
0$~and~$\textbf{a}(x) \in \textbf{W}_2^1(\Omega).$

\textbf{Theorem~4.1.}~~For any right-hand side ~$\textbf{f}(x,t) \in
\textbf{L}_2(Q_t)$~~in equation (4.1) and for any real numbers
~~$\rho > 0 , t > 0,$~~the Navier-Stokes problem-1 has a unique
smooth solution ~$\textbf{u}(x,t) : \textbf{u}(x,t)
\in\textbf{W}_2^{2,1}(Q_t) \cap \textbf{H}_2(Q_t),$ the scalar
function $p(x,t): p_{x_i}(x,t) \in L_2(Q_t)$ satisfying to (4.1)
almost everywhere  on $Q_t,$  and to the initial conditions (4,2).
The following estimate  holds:
$$\|\textbf{u}(x,t)\|_{W^{2,1}_2(Q_t)}~+~\Big\|\frac{\partial p(x,t)}{\partial x_i}\Big\|_{L_2(Q_t)}
\leq~c~\Big(\|\textbf{f}\|_{L_2(Q_t)} +
\|a(x)\|_{W_2^1(\Omega)}\Big)~ \eqno (4.3)$$

$\blacktriangleright$~ In 1941 Hopf  proved that  this problem has a
weak solution $\textbf{u}(x,t)$:  ~
$$\|\textbf{u}(x,t)\|_{L_2(\Omega)} + 2 \rho
\int_0^t\|\textbf{u}_x(x,\tau) \|_{L_2(\Omega)} d \tau <
\|\textbf{a}(x)\|_{W^1_2(\Omega)} + c~ \int_0^t\|\textbf{f}(x,\tau)
\|_{L_2(\Omega)} d \tau
  \eqno(4.4)$$and~~$ \lim_{t\rightarrow
0}\|\textbf{u}(x,t)~-~\textbf{a}(x)\|_{L_2(\Omega)}~=~0.~[3
p.143]$$\blacktriangleleft$

\textbf{The  problem~2.} Find a vector-function
$\textbf{u}^0(x,t)~=~(u_i^0(x,t))_{i=1,2,3} : \Omega \times [0,T]
\rightarrow R^3$   satisfying  the following equation  and the
initial condition:
$$\frac{d \textbf{u}^0(x,t)}{d t} ~-~\rho~ \vartriangle \textbf{u}^0 (x,t)
=~~0;~~~\textbf{u}^0(x,0)~=~\textbf{a}(x),~~\textbf{u}^0(x,t)|_{\partial
\Omega \times [0,T]}~=~~0   \eqno(4,5)$$ $\blacktriangleright$ It
follows by $div~\textbf{a}(x) = 0$ and $\textbf{u}^0(x,t)|_{\partial
\Omega \times [0,T]} = 0$  that: $div ~\textbf{u}^0(x,t) = 0$ for
any t > 0. ~And it follows by $\textbf{a}(x) \in
\textbf{W}_2^1(\Omega)$ that: $|\textbf{u}^0(x,t)| \leq c
\|\textbf{a}(x)\|_{w^1_2(\Omega)}, |\frac{\textbf{u} ^0(x,t)}{d
x_i}| \leq c \|\textbf{a}(x)\|_{w^1_2(\Omega)}.$$\blacktriangleleft$

Then the vector-function: $\textbf{v}(x,t) = \textbf{u}(x,t) -
\textbf{u}^0(x,t)$~ satisfies  the following system of
equations:$$\frac{\partial v_i(x,t)}{\partial t}
~-~\rho~\triangle~v_i(x,t) - \sum_{j=1}^{3}v_j(x,t)~ \frac{\partial
v_i(x,t)}{\partial x_j} - \sum_{j=1}^{3}u_j^0(x,t)~ \frac{\partial
v_i(x,t)}{\partial x_j} - $$$$- \sum_{j=1}^{3}v_j(x,t)~
\frac{\partial u_i^0(x,t)}{\partial x_j} - \sum_{j=1}^{3}u_j^0(x,t)~
\frac{\partial u_i^0(x,t)}{\partial x_j}~+~\frac{\partial
p(x,t)}{\partial x_i}~~=~f_i(x,t)~\eqno(4.6),$$
$$div~\textbf{v}(x,t)~=~\sum_{i=1}^3\frac{\partial v_i(x,t)}{\partial x_i}~=~
0~~,i~=~1,2,3$$and the following initial
conditions:$$\textbf{v}(x,0)~~=~~0,~~~~\textbf{v}(x,t)|_{\partial
\Omega \times [0,T]}~=~~0   \eqno(4,7)$$ Similarly, we introduce the
unknown vector-function ~$\Big(w_i(x,t)\Big)_{i=1,2,3} \in
\textbf{L}_2(Q_t)$~:
$$\frac{\partial v_i(x,t)}{\partial t}
~-~\rho~\triangle~v_i(x,t) ~ - ~\frac{\partial p(x,t)}{\partial
x_i}~~=~w_i(x,t)~\eqno(4.8)$$
$$div~\textbf{v}(x,t)~=~\sum_{i=1}^3\frac{\partial v_i(x,t)}{\partial x_i}~=~
0~~,i~=~1,2,3$$and the following initial
conditions:$$\textbf{v}(x,0)~~=~~0,~~~~\textbf{v}(x,t)|_{\partial
\Omega \times [0,T]}~=~~0   \eqno(4,8')$$In the Proposition 8  we
have proved that for any right-side $\textbf{w}(x,t) \in
\textbf{L}_2(Q_t)$ this problem has a unique solution
$\Big(v_i(x,t)\Big)_{i=1,2,3} \in \textbf{W}_2^{2,1}(Q_t)$~and~
$\|\textbf{v}(x,t)\|_{W_2^{2,1}(Q_t)}~<~c\cdot
\|\textbf{w}(x,t)\|_{L_2(Q_t)}.$ It follows from (4,8) that.
$$v_i(x,t)~=~\int_0^t\int_{\Omega}G(x,t;\xi,\tau)~\Big(w_i(\xi,\tau) + \frac{\partial p(\xi,\tau)}{\partial \xi_i}\Big)
d \xi~d \tau;~~~~i=1,2,3   \eqno(4.9)$$Using (4,8) and this formula,
we rewrite the  Navier-Stokes equation (4,6) as:
$$w_i(x,t) - \sum_{j=1}^{3}G \Big(w_i(\xi,\tau) + \frac{\partial p(\xi,\tau)}{\partial \xi_i}\Big)\cdot
 G_{x_j}\Big(w_i(\xi,\tau) + \frac{\partial p(\xi,\tau)}{\partial \xi_i}\Big)
  - $$$$-\sum_{j=1}^{3}G_{x_j}\Big(w_i(\xi,\tau) + \frac{\partial p(\xi,\tau)}{\partial \xi_i}\Big)
 \cdot u^0_j(x,t) - \sum_{j=1}^{3}G \Big(w_j(\xi,\tau) + \frac{\partial p(\xi,\tau)}{\partial \xi_j}\Big)
\cdot \frac{d u^0_i(x,t)}{d x_j} - \eqno(4,10)$$$$-~
\sum_{j=1}^{3}~u^0_j(x,t)\cdot \frac{d u^0_i(x,t)}{d
x_j}~~~=~~~f_i(x,t)$$ Similarly to the inequality (3,2), using the
estimates (4,4), ~(4,5), it follows  by the equation  (4,10) that:
$$w(t)~~~<~~~f(t) ~~+~~~b~
\Big(\int_0^{t}\frac{w(\tau) + p(\tau)}{(t - \tau)^{\mu}}~d
\tau\Big)^2 ~+~c(T) \cdot \int_0^t \frac{w(\tau) + p(\tau)}{(t -
\tau)^{\mu}}~d \tau
 \eqno(4,11)$$where
$$w(\tau)~ =~\|\textbf{w}(x,\tau)\|_{L_2(\Omega)}~\geq 0;~~~p(\tau)~ =~
\sum_1^3\Big\|\frac{\partial p(x,\tau}{\partial
x_i}\Big\|_{L_2(\Omega)}~\geq 0$$
$$f(t)~ = ~\|\textbf{f}(x,t)\|_{L_2(\Omega)}~+~\sum_{j=1}^{3}\Big\|u^0_j(x,t)\cdot
\frac{d u^0_i(x,t)}{d x_j}\Big\|_{L_2(\Omega)}~ \geq~0
\eqno(4.11')$$
$$c(T)~=~   \|\textbf{a}(x)\|_{L_2(\Omega)} + c~ \int_0^T\|\textbf{f}(x,\tau)
\|_{L_2(\Omega)} d \tau $$ From the inequality $a \cdot b <  a^2 +
b^2$ the following inequality holds:
$$\int_0^t \frac{w(\tau) + p(\tau)}{(t - \tau)^{\mu}}~d \tau ~=~
\int_0^t \frac{w(\tau) + p(\tau)}{(t - \tau)^{\mu/2 }}~ \cdot
\frac{1}{(t - \tau)^{\mu/2}}~d \tau < \int_0^t \frac{(w(\tau) +
p(\tau))^2}{(t - \tau)^{\mu }}d \tau + \frac{t^{1-\mu}}{1-\mu}$$By
this inequality and the inequality  (3,4) in Lemma 3.2 we rewrite
the inequality (4,11) as:
$$w(t) < \Big(f(t) + \frac{t^{1-\mu}}{1-\mu}\Big)~ +~ \Big(b_1  + c(T)\Big)\cdot
\int_0^t \frac{(w(\tau) + p(\tau))^2}{(t - \tau)^{\mu }}~d \tau
$$ Similarly,  based on the proof of Lemma 3.2, we derive the
following inequality:
$$w(t) < \Big(f(t) + \frac{t^{1-\mu}}{1-\mu}\Big) + \Big(b_1 + c(T)\Big)\cdot (1+c)
\int_0^t \frac{w^2(\tau)}{(t - \tau)^{\mu }}~d \tau \eqno(4.12)$$
From the Theorem 3.1,using this inequality, similarly we obtain:~
$$\|\textbf{w}(x,t)\|_{L_2(Q_T)} \leq \sqrt{2}\cdot c(T)\cdot
\|\textbf{f}_u(x,t)\|_{L_2(Q_T)}   \eqno(4,13)$$ where the
vector-function $\textbf{f}_u(x,t)=
\Big(f_i(x,t)\Big)_{i=1,2,3}~+~\Big(\sum_1^3u_j^0(x,t)\frac{d
u_i^0(x,t)}{d x_j}\Big)_{i=1,2,3} $~and c(T) are defined by
(4,11').~~ And we present the basic equation (4,10) as:
$$\textbf{w}(x,t)~-~\Big(K_p + K_1\Big)\ast \Big(E~+~P\Big)\textbf{w}(x,t)
~=~ \textbf{f}(x,t) +~\sum_1^3u^0_j(x,t)\cdot \frac{d
\textbf{u}^0(x,t)}{d x_j} \eqno(4,14)$$ where the
operator~$\Big(K_p\Big)$ is defined by (3,42) , the operator P is
defined by (3,43) and E is an identify operator , i.e. $E
\Big(\textbf{w}(x,t)\Big) = \textbf{w}(x,t).$~The operator
$\Big(K_1\Big)\ast \Big(E~+~P\Big)$~is defined as:
$$K_1 \ast \Big(E+P\Big)\Big(\textbf{w}(x,t)\Big)
= \sum_{j=1}^{3}\Big(\frac{d u^0_i(x,t)}{d x_j}\cdot G +
u^0_j(x,t)\cdot G_{x_j}\Big)\ast
\Big(E+P\Big)\Big(\textbf{w}(x,t)\Big) \eqno(4,15)$$ As the proof of
Theorem 2.1(p.19) and definitions of the operators $K_p, P , K_1$ by
the formulas (3,42), (3,43) and (4,15) it is proves that the
following operators: ~~~$\Big(K_p \Big)\ast \Big(E~+~P\Big)$;

$\Big(K_1 \Big)\ast \Big(E~+~P\Big)$~are compact on the vector-space
$\textbf{L}_2(Q_t).$ ~~Similarly, using  Lerau-Schauder's theorem
and the estimate (4.13), the existence and smoothness of solution to
the Navier-Stokes problem with $u(x,0) = a(x)\neq 0$ is
proved.$\blacktriangleleft$~ ~~Thank you for your attention. \vskip
5mm

\section*{Declarations}

\subsection*{Competing interests} The authors declare that they have
no competing interests.

\subsection*{Funding} There are no funding sources for this
manuscript.

\subsection*{Acknowledgments} Not applicable.

\newpage

\begin{center}  References
\end{center}

[1].~ Leray J.~: Sur le Mouvement d'un Liquide Emplissent l'Espace ,
Acta.Math.

~~~~~~51(1934).

[2].~Temam P.~: Navier-Stokes equations , theory and numerical
analysis ,Oxford ,(1979).

[3].~ Ladyzhenskaya O.~:The Mathematical Theory of Viscous
Incompressible Flows,

~~~~~~Gordon and Breach , (1969) .

[4].~Lions J.L.~:Quelques methodes de resolution des problemes aux
limites non lineaires,~

~~~~~~Paris , (1969).

[5].~Ladyzhenskaya O.,Solonnikov B., Uralzeva N.~: The linear and
quasi-linear equations

~~~~~~of parabolic type , Moscow , (1967).

[6].~ Ladyzhenskaya O.~:~Problems for Navier - Stokes equation
~,~~Successes of

~~~~~~math sciences ~/ ~~2003~,~v.58~,~2(350) ~,~March -April .

[7].~ Tricomi F.G.~:~Integral equations.  , Interscience publishers
, INC.,New York (1957).

[8].~ Samko S.,Kilbas A.,Marichev O.~: Fractional Integrals and
Derivative . ~Theory

~~~~~~and Applications. , Cordon and Breach Sci.Publishers , (1993).

[9]. Friedman Avner .~:Partial differential equations of parabolic
type., ~

~~~~~~Prentice-Hall , (1964).

[10].Kamke E:  Gewohnliche Differential  Gleichungen. , Leipzig ,
(1959).

[11].Sobolev S.L.Partial differential equations of Mathematical
Physics.,Moskow,(1966).

[12]. Bazarbekov A. "Estimates of solutions to the linear
Navier-Stokes equation, 5 pages

~~~~~~in Latex". ~This manuscript has not been  published.

\vskip35mm

Address:050026 ,Muratbaeva St. 94-5,

Almaty , Republic Kazakhstan.\label{BAZARBFIN}

\end{document}